\documentclass{siamltex}

\usepackage{amsmath,amsfonts,amssymb}
\usepackage{subfig}
\usepackage[rightcaption]{sidecap}
\usepackage{enumerate}
\usepackage{url}
\usepackage{graphics,graphicx}
\usepackage{rotating}
\usepackage{placeins}
\usepackage{algorithm,algorithmic}
\usepackage{array,longtable}
\usepackage[table]{xcolor}
\usepackage{color}
\usepackage{todonotes}

\newcommand{\btheta}{{\boldsymbol \vartheta}}
\newcommand{\e}[1]{\ensuremath{{\times 10^{#1}}}}

\usepackage[pdftex]{hyperref}
\hypersetup{colorlinks,citecolor=black,
filecolor=black,linkcolor=black,urlcolor=black}

\title{Universal Meshes: A new paradigm for computing with
  nonconforming triangulations}

\author{Ramsharan Rangarajan and Adri\'{a}n J. Lew}
\date{}
\begin{document}

\maketitle

\begin{abstract} 
  We describe a method for discretizing planar $C^2$-regular domains
  immersed in non-conforming triangulations. The method consists in
  constructing mappings from triangles in a background mesh to
  curvilinear ones that conform exactly to the immersed
  domain. Constructing such a map relies on a novel way of
  parameterizing the immersed boundary over a collection of nearby
  edges with its closest point projection. By interpolating the
  mappings to curvilinear triangles at select points, we recover
  isoparametric mappings for the immersed domain defined over the
  background mesh. Indeed, interpolating the constructed mappings just
  at the vertices of the background mesh yields a fast meshing
  algorithm that involves only perturbing a few vertices near the
  boundary.
    
  For the discretization of a curved domain to be robust, we have to
  impose restrictions on the background mesh. Conversely, these
  restrictions define a family of domains that can be discretized with
  a given background mesh. We then say that the background mesh is a
  {\it universal mesh} for such a family of domains. The notion of
  universal meshes is particularly useful in free/moving boundary
  problems because the same background mesh can serve as the universal
  mesh for the evolving domain for time intervals that are independent
  of the time step. Hence it facilitates a framework for finite
  element calculations over evolving domains while using a fixed
  background mesh. Furthermore, since the evolving geometry can be
  approximated with any desired order, numerical solutions can be
  computed with high-order accuracy. We demonstrate these ideas with
  various numerical examples.
\end{abstract}


\section{Introduction}
\label{sec:introduction}
Finite element methods commonly handle evolving domains in one of two
ways --- either the changing domain is remeshed at each
instant/update, or it is immersed in a background mesh and
approximated within it. We introduce a novel approach here that
inherits the conceptual simplicity of the former and the computational
efficiency of the latter. We describe a method for discretizing
sufficiently smooth planar domains \emph{using} a given background
triangulation, provided some conditions are met.  The method consists
in constructing mappings from triangles in a background mesh to
curvilinear ones that conform exactly to the immersed domain.  As an
example, consider simulating a problem in which rigid blades
physically mix fluid in a closed container. As the blades rotate, the
region of the container occupied by the fluid changes. We can now
discretize the evolving fluid domain by merely perturbing vertices and
edges of the \emph{same} background mesh. Precisely because the same
background mesh is utilized for all positions of the propeller, we
term it a \emph{universal mesh} for the fluid. Since connectivities of
triangles remain unaltered and no new vertices are introduced, sparse
patterns of data structures involved in the problem can also be
retained.

{\sc The basic idea}: The key idea in constructing conforming
discretizations for an immersed domain consists of a two-step
procedure to perturb edges and vertices in its background mesh. We
identify triangles in the background mesh with at least one vertex
inside the curved domain. Edges belonging to this collection having
both vertices outside the domain are mapped onto the boundary with its
closest point projection.  Then we relax away from the boundary a few
vertices that lie inside the domain and close to the boundary. For a
certain class of background meshes, these steps enable us to construct
a homeomorphism between the union of selected triangles in the
background mesh and the the domain.  In other words, this construction
yields a conforming curvilinear discretization for the immersed
domain.

There are no conformity requirements on the background mesh; for
instance, none of its vertices need to lie on the boundary. The
resulting algorithmic advantages are significant, especially for
problems with evolving domains. Such problems are ubiquitous,
including ones with interaction between fluids and solids, problems
with free boundaries and moving interfaces, domains with propagating
cracks and problems with phase transformations. Various numerical
schemes have been proposed for such applications, see references
\cite{bazilevs2008isogeometric, donea2004arbitrary, osher2003level,
  peskin2002immersed, wagner2001extended} for a representative few.
The spirit of this article, as evidenced by the examples presented, is
to immerse such evolving geometries in a universal mesh and update its
spatial discretization as necessary.

Numerical methods that adopt nonconforming meshes have been formulated
in various ways. For instance, by re-triangulating elements near the
boundary and with cut/trimmed cells \cite{moumnassi2011finite,
  saksono2007adaptive}; by treating immersed boundaries and interfaces
via constraints using penalty \cite{angot1999penalization}, Lagrange
multipliers \cite{burman2010fictitious} and Nitsche's method
\cite{hansbo2002unfitted}; or by enriching the space of solutions near
the boundary as done in extended finite element methods and
discontinuous Galerkin-based methods \cite{wagner2001extended,
  lew2008discontinuous}.  One of the challenges in these methods is
achieving optimal accuracy with high-order interpolations. Almost
without exception, these methods resort to a polygonal approximation
for the immersed domain. Such approximations suffice in low order
methods, in which the solution is approximated by piecewise constant
or affine functions \cite{lew2008discontinuous,
  rangarajan2009discontinuous}. To construct high-order methods, it is
imperative to approximate the immersed geometry sufficiently well over
the background mesh.

With conforming curvilinear discretizations for immersed domains, we can 
construct curved finite elements with optimal convergence rates given
only nonconforming background meshes (\S\ref{sec:curved-elements}).
Alternately, by interpolating the mappings to curvilinear triangles at
select points, we recover isoparametric mappings defined over the
background mesh, as discussed in \S\ref{sec:isoparametric}.  Rather
than discretizing a curved domain exactly, isoparametric mappings
provide a systematic way of constructing \emph{sufficiently accurate}
approximations in which curved edges interpolate the boundary at
select points/nodes. Both exactly conforming and isoparametric
elements enable high-order convergence rates (see
\S\ref{subsec:cvg}). The former type of elements can be particularly
advantageous in problems sensitive to boundary conditions. One such
problem is that of a simply supported circular plate in bending
\cite{babuska1990plate}. In \S\ref{subsec:plate}, we show that the
accuracy of its numerical solution can be sensitive to how well curved
boundaries are represented.

The construction for curved domains given here are inspired by the
well known mappings proposed in references
\cite{gordon1973transfinite, lenoir1986optimal,
  mansfield1978approximation, scott1973finite, zlamal1973curved}. The
constructions in these articles assume as a point of departure, (i) a
conforming mesh for the curved domain, and (ii) a parametrization for
the curved boundary over edges that interpolate it. By admitting
nonconforming meshes and not relying on a specific representation for
the boundary, we relax both assumptions. Hence we generalize these
known constructions to a larger class of meshes and amenable to
general boundary representations. For example, the boundary can be
given parametrically as a collection of splines, or implicitly as the
zero level set of a function. We only require a way to compute the
closest point projection sufficiently close to the boundary.

In the special case when the mappings to curvilinear triangles are
interpolated just at vertices in the background mesh, we get a
conforming mesh for the immersed domain. The resulting ``meshing
algorithm'' is very fast, since it only involves perturbing
vertices. For simplicity, we present the meshing algorithm first, in
\S\ref{sec:meshing}, by defining affine mappings that perturb vertices
of the background mesh near the immersed boundary. In contrast,
conforming curvilinear discretizations, discussed in
\S\ref{sec:curved-elements}, are constructed by mapping certain
\emph{edges}, rather than just vertices, onto the boundary.

It is common knowledge that perturbing vertices in a background mesh
\emph{can} yield a conforming mesh for an immersed domain. The
challenge, though, is making such algorithms robust. An algorithm
specific to the case of domains immersed in rectangular grids is
described in \cite{borgers1990triangulation}.  The closest point
projection has also been used to locally modify Cartesian grids near
the boundary, as done in \cite{gonzalez2006inverse,
 sanches2011immersed} although the authors do not consider when their
approach could fail.

In \S\ref{subsec:assumptions}, we list sufficient conditions for our
meshing algorithm to be robust. We assume that the curved domain is
$C^2$-regular, and we have to \emph{restrict the class of background
  meshes} for a given domain. In particular, we require that the
background mesh be sufficiently refined and that certain angles in
triangles near the boundary be acute. With these assumptions, we
analyzed the restriction of the closest point projection to the edges
whose vertices are projected onto the immersed boundary in
\cite{rangarajan2011analysis}. We proved that this mapping is in fact
a homeomorphism onto the boundary. As we discuss in
\S\ref{sec:rationale}, such a result and a possibly smaller mesh size
ensures that moving vertices of the background mesh in the way we do
in the meshing algorithm will not result in degenerate, inverted, or
overlapping triangles, and in general, avoids tangled meshes. For
instance, a refined background mesh of equilateral triangles is
guaranteed to mesh a smooth domain immersed in it. We cannot however
make the same claim with a background mesh of right angled triangles,
because such a mesh may not satisfy the condition on angles. The
conditions required for the success of the meshing algorithm along
with possibly smaller mesh size near the boundary also ensure that the
mappings to conforming and isoparametric curved elements are well
defined as well.

\section{From background meshes to conforming meshes}
\label{sec:meshing}
We begin by illustrating the steps to determine a conforming mesh for
a planar curved domain $\Omega$ immersed in a background triangulation
${\cal T}_h$, where $h$ denotes the mesh size. We assume that $\Omega$
is an open set and denote $\Omega^c = {\mathbb R}^2\setminus
\Omega$. By $\Omega$ being \emph{immersed} in ${\cal T}_h$, we mean
that the set triangulated by ${\cal T}_h$ contains
$\overline{\Omega}$. To ensure that the resulting mesh for $\Omega$ is
valid, we require certain assumptions on $\Omega$ and ${\cal
  T}_h$. They are stated in \S\ref{subsec:assumptions} and discussed
in \S\ref{sec:rationale}.
\begin{center}
  \begin{longtable}{m{1.7in} m{2.9in}}
    \caption{Steps in the meshing algorithm} \label{algo}
    \\ \hline \\
    
    \textbf{Step 1:} Identify vertices in the background mesh ${\cal
      T}_h$ that lie in $\Omega$ (square markers) and in $\Omega^c$
    (circular markers) respectively. Omit triangles with no vertices
    in $\Omega$.
    
    &\includegraphics[scale=0.5]{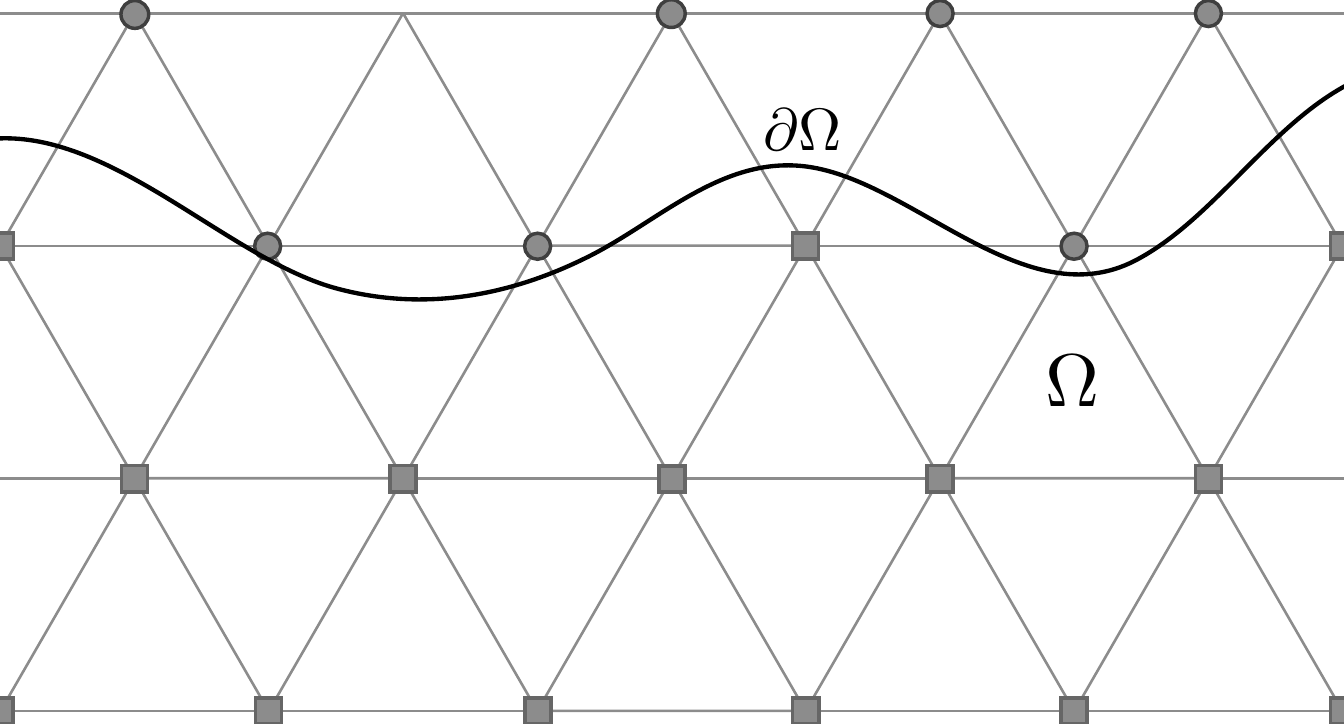} \\[5pt]
    
    \textbf{Step 2:} Identify \emph{positively cut triangles} in
    ${\cal T}_h$ --- triangles with precisely one vertex in
    $\Omega$. These are shaded in gray. In each such triangle, check
    that the angle at its vertex in $\Omega^c$ closest to
    $\partial\Omega$ is smaller than $90^\circ$. These angles are
    labeled $\vartheta$ in the adjacent figure.
    
    & \includegraphics[scale=0.5]{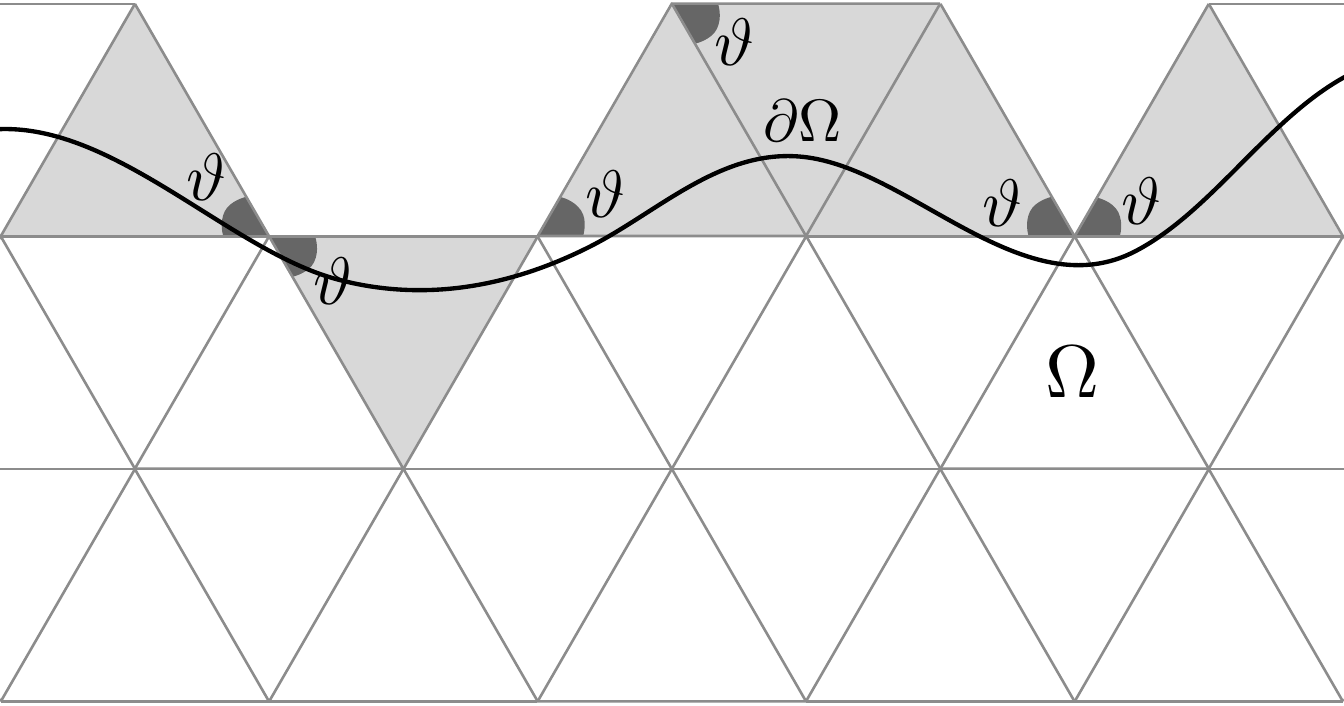} \\[5pt]
    
    \textbf{Step 3:} A \emph{positive edge} is the edge of a
    positively cut triangle joining its vertices in $\Omega^c$. These
    edges are shown in black. Map vertices of positive edges to their
    closest point on the boundary $\partial\Omega$.
    
    & \includegraphics[scale=0.5]{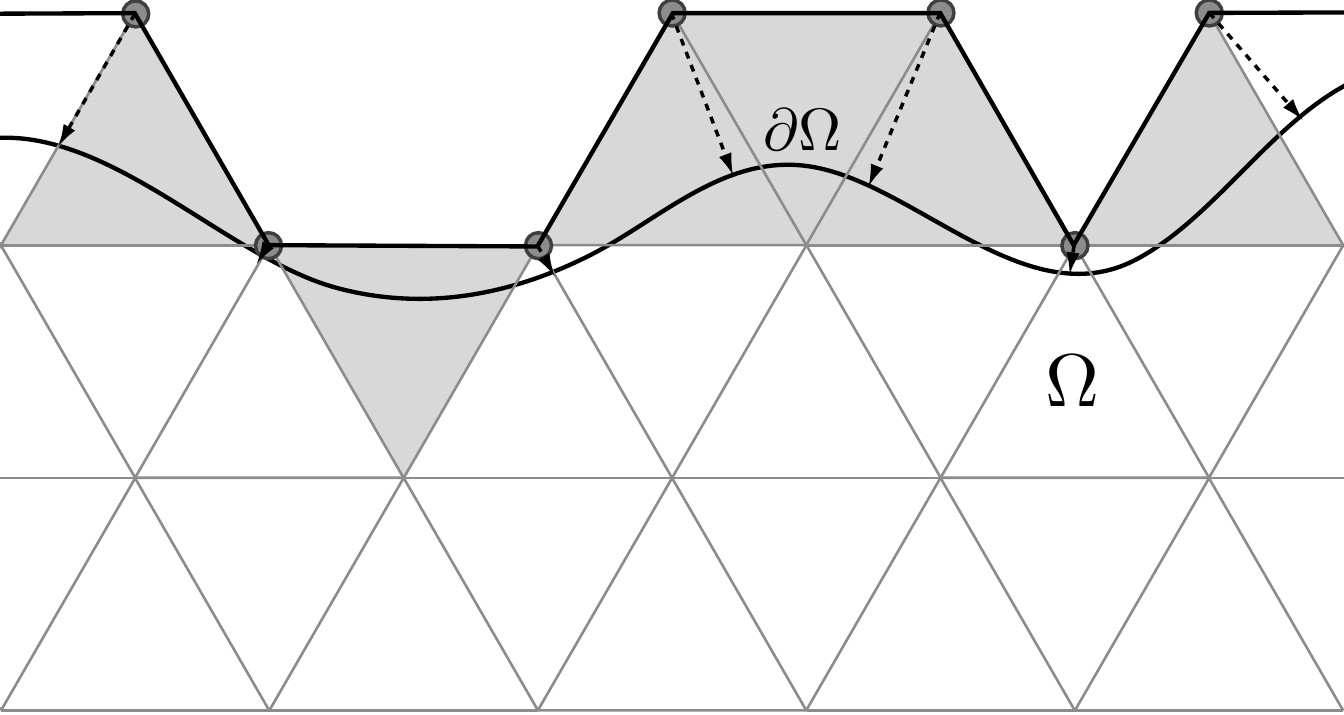} \\[5pt]

    \textbf{Step 4:} Identify vertices in $\Omega$ that lie close to
    $\partial\Omega$, such as the ones shown by triangular
    markers. Perturb these vertices away from $\partial\Omega$ using
    the mapping $\mathfrak{p}_h$ in \eqref{eq:pert}.
    
    & \includegraphics[scale=0.5]{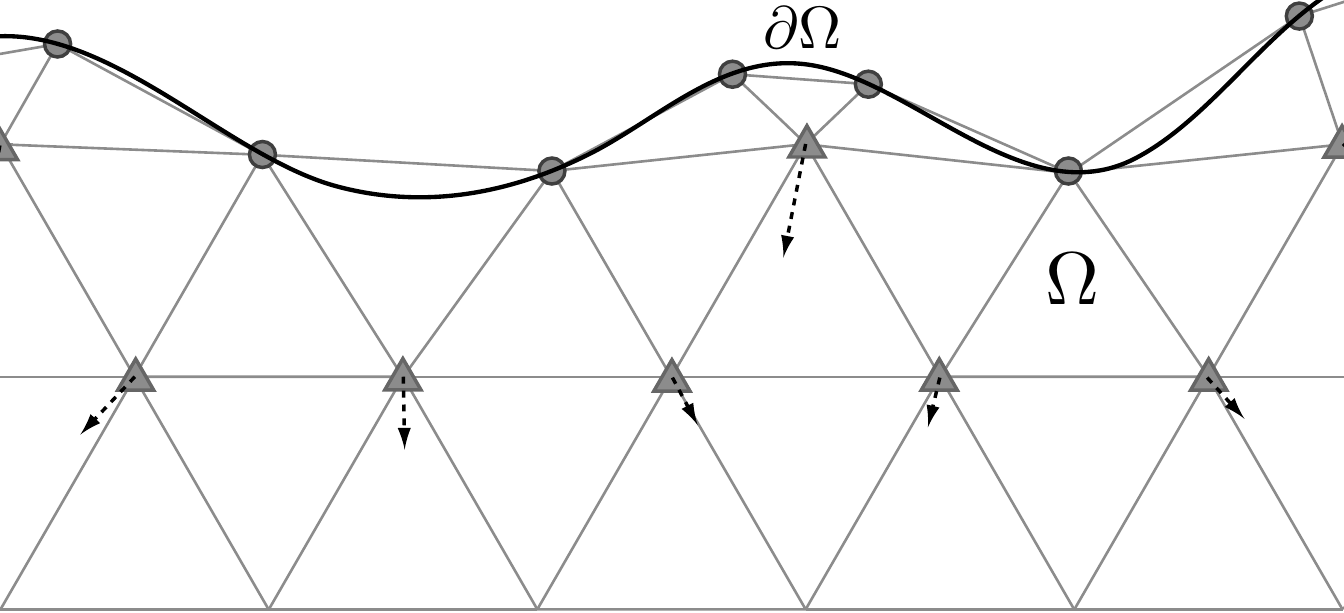} \\[5pt]
    {\sc \textbf{ Final Mesh}} &
    \includegraphics[scale=0.5]{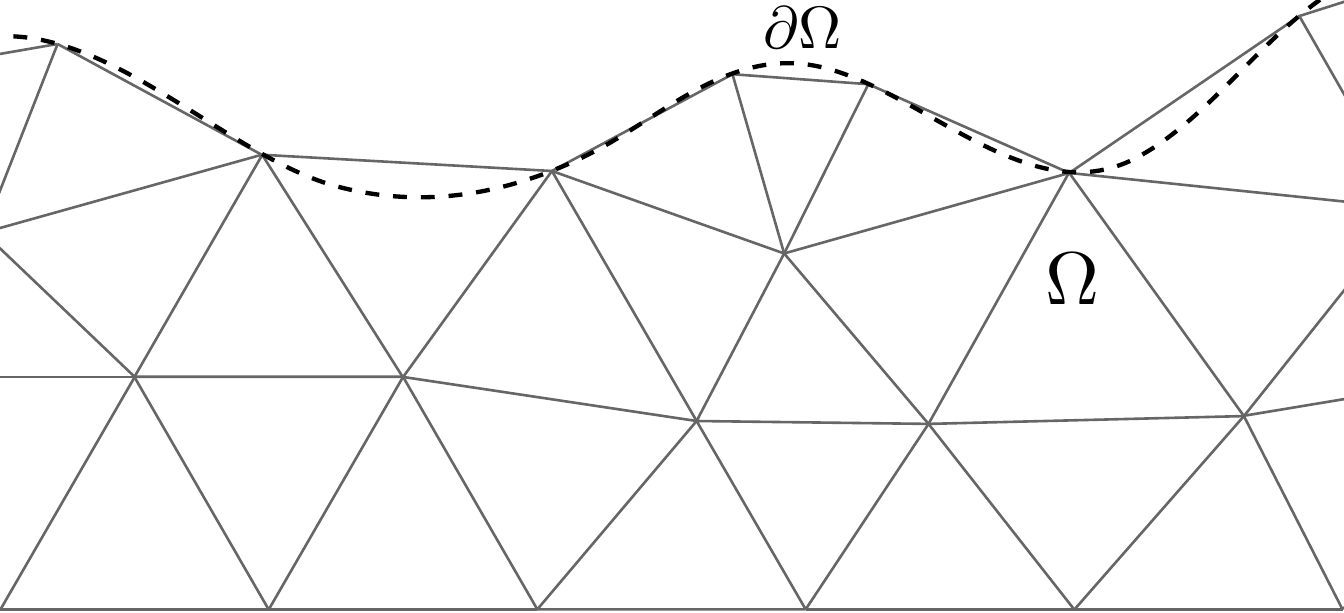} \\[7pt]
    \hline
  \end{longtable}
\end{center}

The mapping for relaxing vertices away from the boundary is given by
\begin{align}
  \mathfrak{p}_h(x) = \begin{cases}
    x - \eta h\left(1+\frac{\phi(x)}{r}\right) \nabla\phi(x) ~ ~ \text{if}~ -r <\phi(x)<0, \\
    x ~ ~ \text{otherwise},
    \end{cases}
    \label{eq:pert}
  \end{align}
  where $\phi$ is the signed distance function to $\partial\Omega$
  defined below in \S\ref{subsec:details}, $\eta\in (0,1)$ and $r$
  equals a few multiples of the mesh size $h$. We discuss the choice
  of $\eta$ and $r $ in \S\ref{subsec:details} as well.

  The meshing algorithm in Table \ref{algo} is succinctly summarized as
  a piecewise affine mapping over triangles in ${\cal T}_h$. To this
  end,  let $\pi:{\mathbb
  R}^2\rightarrow \partial\Omega$ denote the closest point projection onto
$\partial\Omega$, and identify the collections of triangles
  \begin{align}
    {\cal T}_h^i &= \{K\in {\cal T}_h\,:\,\phi\geq 0~\text{at
      precisely}~i~\text{vertices of}~K\}
\end{align}
for $ i=0,1,2,3$. Consider a triangle $K\in {\cal T}_h^{0,1,2}$ with
vertices $\{u,v,w\}$ ordered such that $\phi(u)\geq \phi(v)\geq
\phi(w)$. Denote the barycentric coordinates of $x\in \overline{K}$ by
$(\lambda_u,\lambda_v,\lambda_w)$ so that
$x=\lambda_u\,u+\lambda_v\,v+\lambda_w\,w$ and $\lambda_u+\lambda_v+\lambda_w=1$. The algorithm in Table \ref{algo} maps $x\mapsto
M_K^h(x)$ defined as
\begin{align}
  M_K^h(x) &= 
  \begin{cases}
    \lambda_u\mathfrak{p}_h(u)+\lambda_v\mathfrak{p}_h(v)+\lambda_w\mathfrak{p}_h(w) ~ &\text{if}~K\in {\cal T}_h^0, \\
    \lambda_u\pi(u)+\lambda_v\mathfrak{p}_h(v)+\lambda_w\mathfrak{p}_h(w) ~ &\text{if}~K\in {\cal T}_h^1, \\
    \lambda_u\pi(u)+\lambda_v\pi(v)+\lambda_w\mathfrak{p}_h(w) ~ &\text{if}~K\in {\cal T}_h^2. 
  \end{cases}
  \label{eq:Mk}
\end{align}

To refer to the angles checked in step 2 in Table \ref{algo}, we
introduce the terms \emph{proximal vertices} and \emph{conditioning
  angles}. The proximal vertex of a positively cut triangle is the
vertex of its positive edge closer to $\partial\Omega$. If both
vertices of the positive edge are equidistant from $\partial\Omega$,
the one containing the smaller interior angle is chosen as the
proximal vertex. If the angles are equal as well, the choice is
arbitrary. The conditioning angle of a positively cut triangle is the
interior angle at its proximal vertex. Hence in step 2 of Table
\ref{algo}, we require that each positively cut triangle have an acute
conditioning angle.

\subsection{Sufficient conditions for a valid mesh}
\label{subsec:assumptions}
We require a few assumptions to guarantee that the meshing algorithm
in Table \ref{algo} yields a valid mesh for the immersed domain
$\Omega$. By a valid mesh we mean a triangulation of a polygon
$\Omega_h$ made of triangles with diameters smaller or equal to $h$,
such that the vertices of $\Omega_h$ lie on $\partial \Omega$, and
$\Omega_h$ approximates $\Omega$ as $h\searrow 0$\footnote{At least we
  need to have $|\Omega_h\setminus \Omega\cup
  \Omega\setminus\Omega_h|\to 0$ and
  $\text{distance}(\partial\Omega_h,\partial\Omega)\to 0$ as
  $h\searrow 0$. }.  We require that
\begin{enumerate}[(a)]
\item the domain $\Omega$ be $C^2$-regular,
\item that $\Omega$ be immersed in ${\cal T}_h$, i.e.,
  $\overline{\Omega} \subset \cup_{K\in{\cal T}_h}\overline{K}$,
\item the conditioning angle in each positively cut triangle in ${\cal
    T}_h$ be strictly smaller than $90^\circ$, and
\item the triangulation ${\cal T}_h$ be sufficiently refined in the
  vicinity of $\partial\Omega$.
\end{enumerate}
A precise definition of $C^k$-regular domains is given in
\cite{rangarajan2011analysis}. For our purposes, it suffices to note
that $\Omega$ is $C^2$-regular if the signed distance function to
$\partial\Omega$ is $C^2$ in a neighborhood of
$\partial\Omega$. Assumption (d) requires that the mesh size be small
near $\partial\Omega$. By this we mean that if triangle $K\in{\cal
  T}_h$ lies near $\partial\Omega$, then its diameter $h_K$ should be
smaller than a value that depends on the local curvature and feature
size of $\partial\Omega$, among others.  Explicit estimates for the
required mesh size near $\partial\Omega$ are essential to automate the
algorithm. However, we do not provide such estimates here. We do
briefly mention how to construct adaptively refined background meshes
satisfying the acute conditioning angle requirement (c) in
\S\ref{subsec:details}. We discuss the rationale behind these
assumptions later in \S\ref{sec:rationale}.

\subsection{An illustrative example}
\begin{figure}
  \centering \subfloat[\label{fig:bunny-a}A curved domain is immersed
  in a background mesh of equilateral triangles. The boundary of the
  domain shown is composed of $41$ cubic splines. Its is therefore
  $C^2$-regular. Since all interior angles of triangles in the
  background mesh equal $60^\circ$, the conditioning angle also equals
  $60^\circ$. Triangles with at least one vertex inside the domain
  (i.e., triangles in ${\cal T}_h^{0,1,2}$) are shaded in gray.  These
  triangles are mapped to a conforming mesh for the immersed domain.]
  {\includegraphics[scale=0.45]{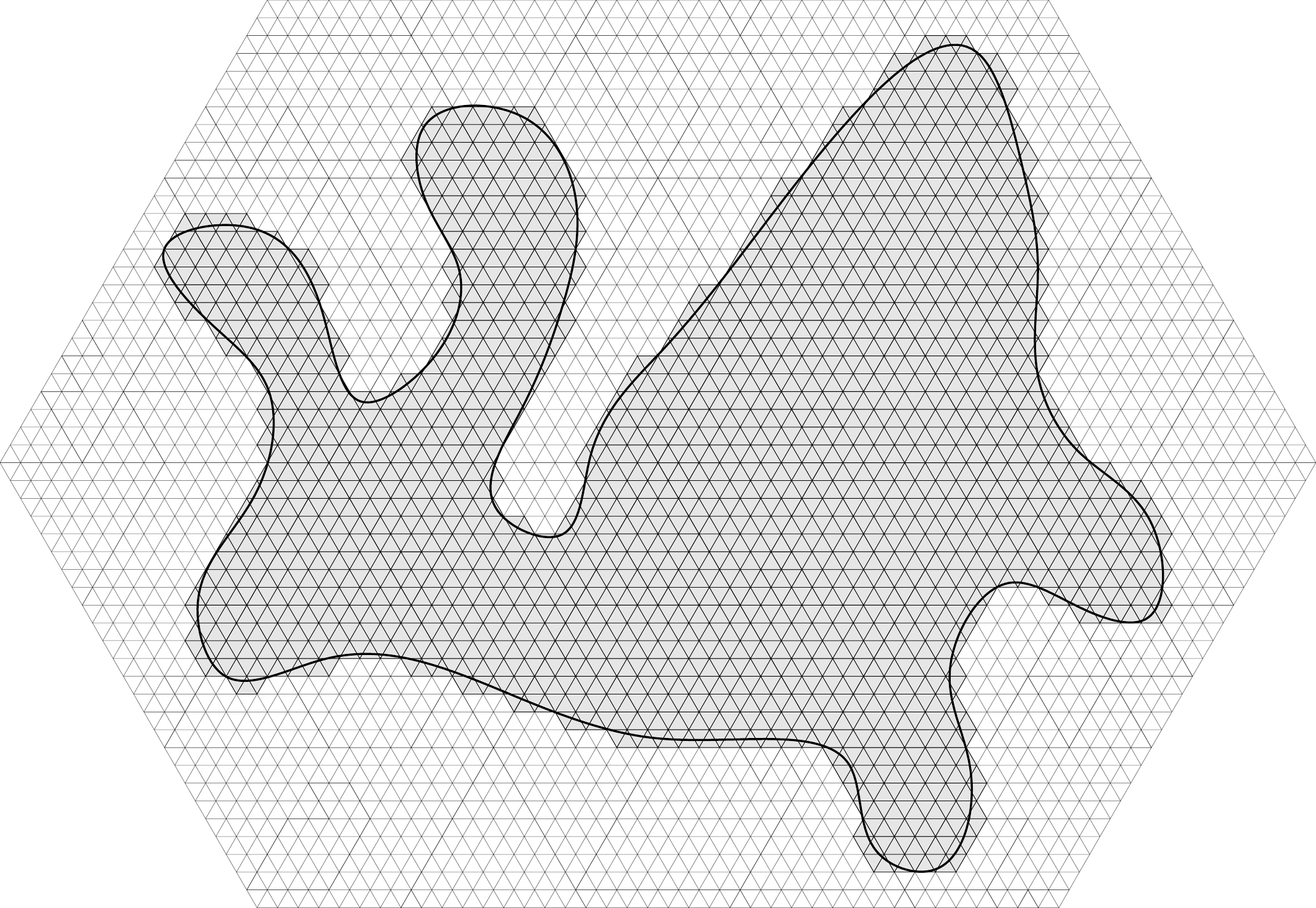}}
  \\
  \subfloat[\label{fig:bunny-b}There are $2382$ triangles in the
  collection ${\cal T}_h^{0,1,2}$. Of these, $234$ triangles belong to
  ${\cal T}_h^2$ (shaded in black), $228$ triangles to ${\cal T}_h^1$
  (shaded in dark gray) and the remaining $1920$ triangles to ${\cal
    T}_h^0$. The $756$ triangles left unshaded in the figure are the
  ones that are retained unaltered from the background mesh. The
  remaining ones have their vertices either snapped onto the boundary
  or relaxed away from the boundary. ]
  {\includegraphics[scale=0.4]{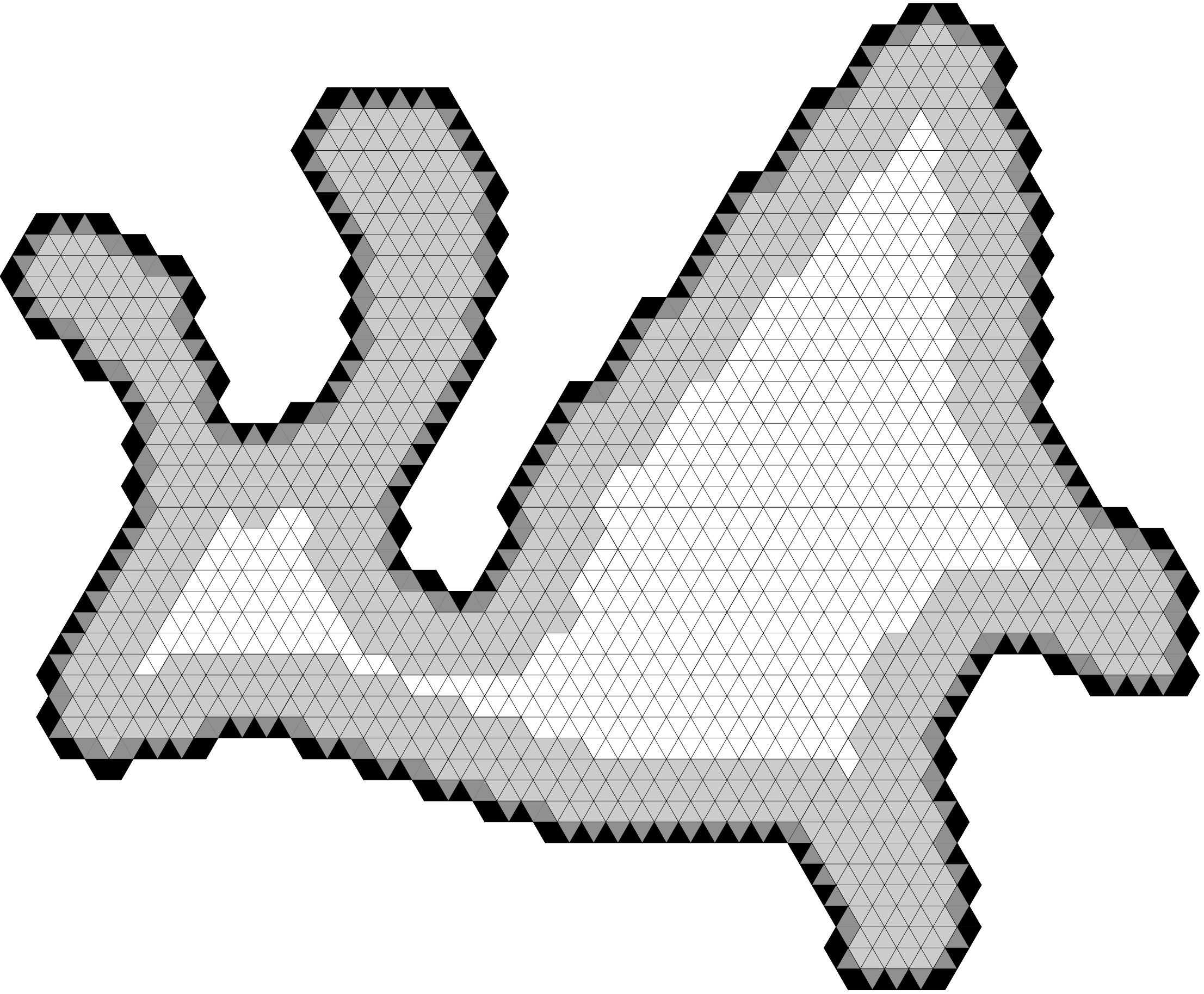}}
  \caption{Example illustrating the algorithm in Table \ref{algo} to
    mesh an immersed domain by perturbing vertices in a background
    mesh.}
  \label{fig:bunny}
\end{figure}

\begin{figure}
  \centering
  \subfloat[\label{fig:bunnymesh-1}]{\includegraphics[scale=0.7]{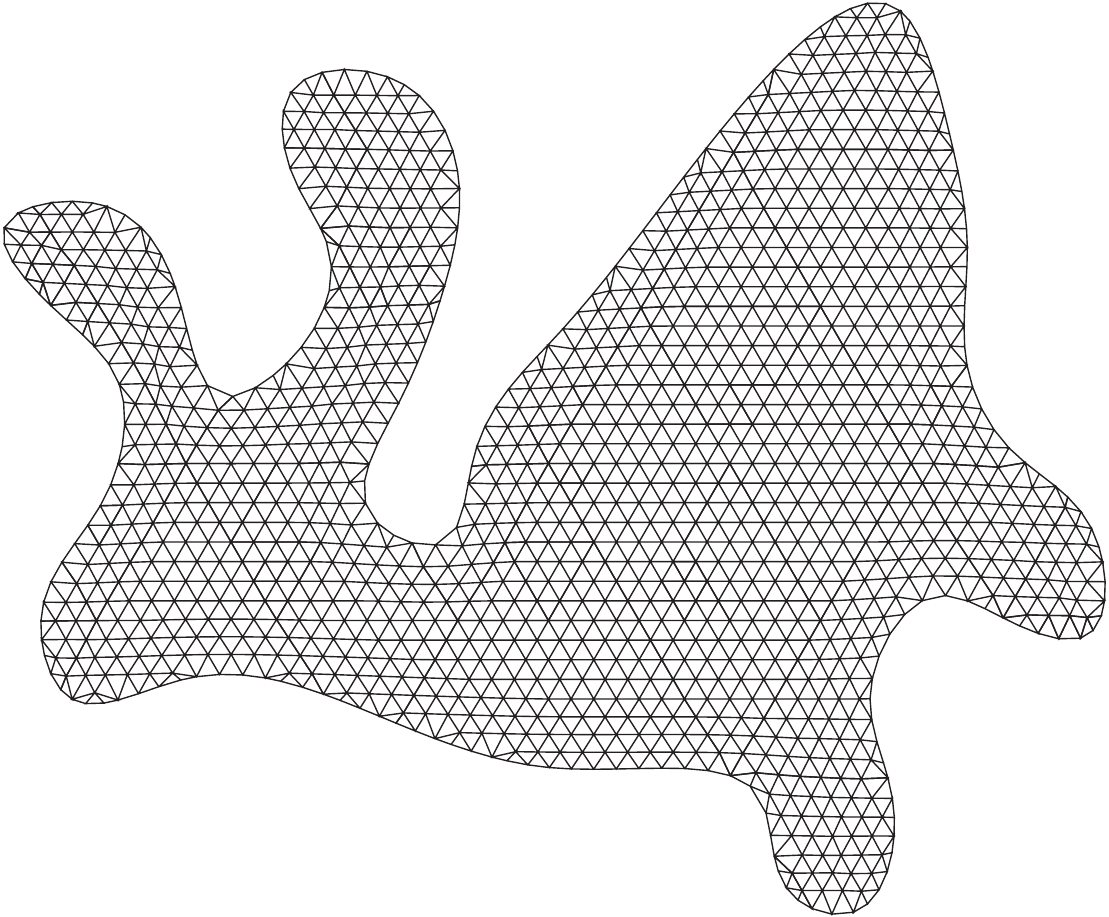}} 
  \hfill
  \subfloat[\label{fig:bunnymesh-2}]{\includegraphics[scale=0.9]{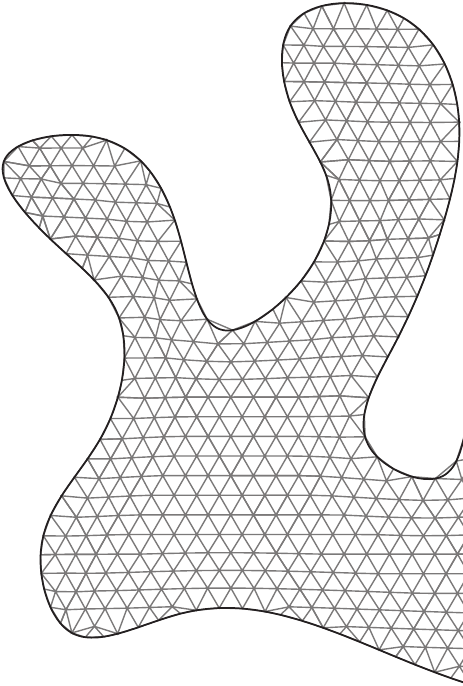}}
  \caption{The conforming mesh determined by the meshing algorithm for
    the domain and background mesh shown in Fig. \ref{fig:bunny}. A
    closer view is shown in (b).}
  \label{fig:bunnymesh}
\end{figure}
An example that uses the meshing algorithm is shown in
Fig. \ref{fig:bunny}. The curved domain to be meshed is the one shown
in Fig. \ref{fig:bunny-a}. It is $C^2$-regular because its boundary is
a collection of cubic splines. It is immersed in a background mesh of
equilateral triangles. Hence the conditioning angle equals $60^\circ$
and the check in step 2 in Table \ref{algo} (also assumption (c) in
\S\ref{subsec:assumptions}) is trivially satisfied. In
Fig. \ref{fig:bunny-a}, the $2382$ triangles in ${\cal T}_h^{0,1,2}$
are shaded in gray. Triangles in this collection are mapped to a
conforming mesh for the immersed domain as shown in
Fig. \ref{fig:bunnymesh}. Of these triangles, $756$ remain unaltered
from the background mesh. At least one vertex in the remaining
triangles is perturbed.  Note that with a more refined background
mesh, a larger fraction of the triangles in ${\cal T}_h^{0,1,2}$ will
remain unaltered from the background mesh. For instance, when the
background mesh in Fig. \ref{fig:bunny} is refined once by
subdivision, $5552$ of the $9178$ triangles in the resulting
conforming mesh remain unaltered (i.e., remain equilateral triangles).

In Table \ref{histo}, we inspect the quality of triangles in the mesh
in Fig. \ref{fig:bunnymesh}. We use the ratio of the circumradius to
the inradius as a metric for the quality of triangles. The best
possible value of this ratio is $2$, which is attained in equilateral
triangles. The table lists the number of triangles in the final mesh
with quality in a given range of the metric.  The minimum and maximum
angles in the mesh were $20.6^\circ$ and $129.6^\circ$. Table
\ref{histo} also reports the quality of the mesh determined by the
algorithm upon refining the background mesh in Fig. \ref{fig:bunny-a}
by subdividing each triangle into four self-similar ones. Extreme
angles in the resulting mesh for the curved domain were $18.4^\circ$
and $139.7^\circ$ respectively.

\begin{table}
\centering
\caption{Quality of the mesh determined by the meshing algorithm for
  the domain in Fig. \ref{fig:bunny-a}. The metric used for the quality
  of a triangle is the ratio of its circumradius to
  the inradius. We
  only inspect triangles in the final mesh that have been perturbed
  with respect to the background mesh by the meshing algorithm. The
  remaining triangles remain equilateral. The
  column titled `coarse mesh' lists the number of triangles in the
  mesh in Fig. \ref{fig:bunnymesh} that have quality in the range
  specified in the first column. The quality ranges from $2.0$ to
  $5.8$. The column titled `refined mesh'
  lists corresponding values for the mesh determined using a
  self-similar refinement of the background mesh shown in
  Fig. \ref{fig:bunny-a}. In this case, the quality ranges from $2.0$ to
  $8.8$. }
\label{histo}
\begin{tabular}{|c|c|c|}
  \hline
  Range of metric & coarse mesh & refined mesh\\  \hline
  $2.0-2.4$ & $1530$ & $3441$ \\ \hline
  $2.4-2.8$ & $45$   & $75$\\ \hline
  $2.8-3.2$ & $16$   & $42$\\ \hline
  $3.2-3.6$ & $9$    & $22$\\ \hline
  $3.6-4.0$ & $6$    & $11$\\ \hline
 $4.0-4.4$ & $7$     & $7$\\ \hline
 $4.4-4.8$ & $5$     & $10$\\ \hline
 $4.8-5.2$ & $4$     & $3$\\ \hline
 $5.2-5.6$ & $2$     & $2$\\ \hline
 $5.6-6.0$ & $2$     & $3$\\ \hline
 $6.0-6.4$ & $0$     & $2$\\ \hline
 $6.4-6.8$ & $0$     & $4$\\ \hline
 $>6.8$     & $0$    & $4$\\ \hline
\end{tabular}
\end{table} 

\subsection{Background meshes as \emph{Universal meshes}}
\label{subsec:universal}
An important advantage of admitting nonconforming background meshes is
in problems with evolving domains. For then it is possible, at least
in principle, to use the same background mesh to triangulate a
changing domain. If not for the entire duration of interest, at least
for reasonably large changes in the immersed geometry. This motivates
the notion of \emph{universal meshes}.

\begin{figure}
  \centering
  \subfloat[\label{fig:bg-splines}]{\includegraphics[scale=0.75]{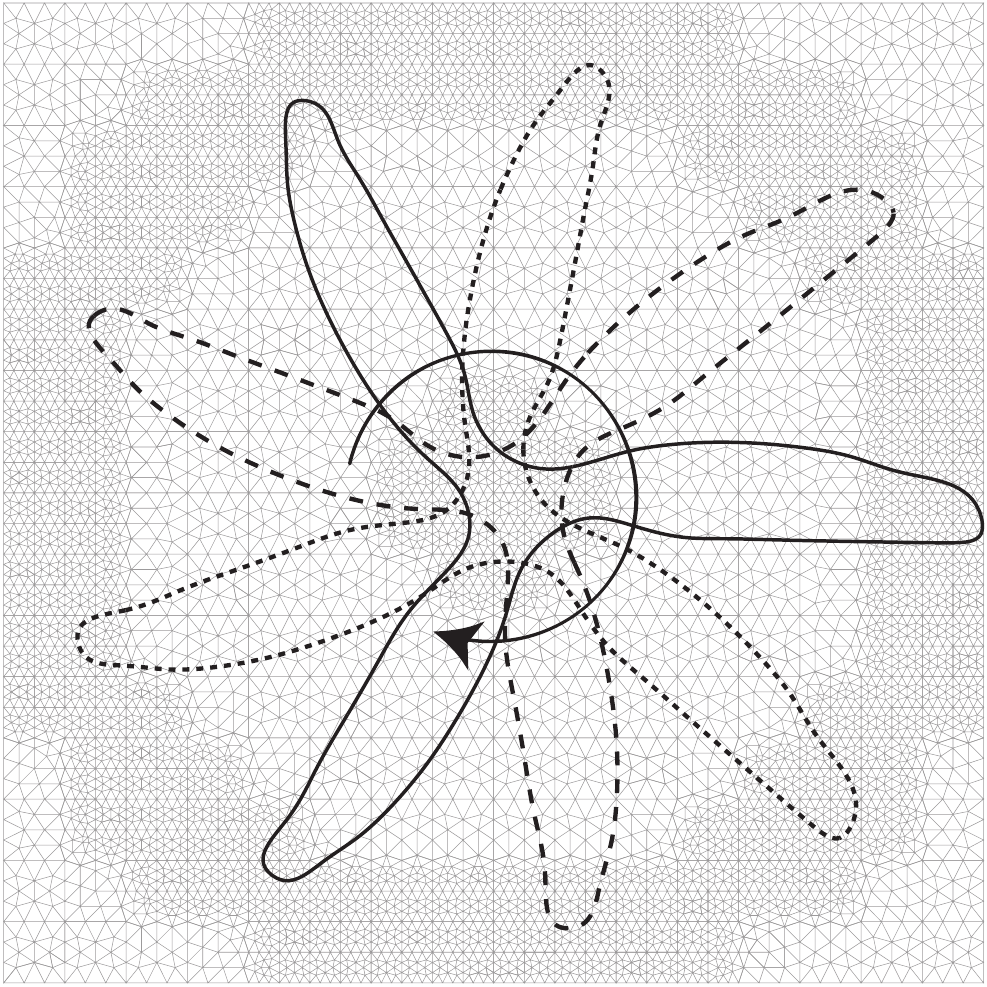}}\\
  \subfloat[\label{fig:propeller-meshes}]{\includegraphics[scale=0.6]{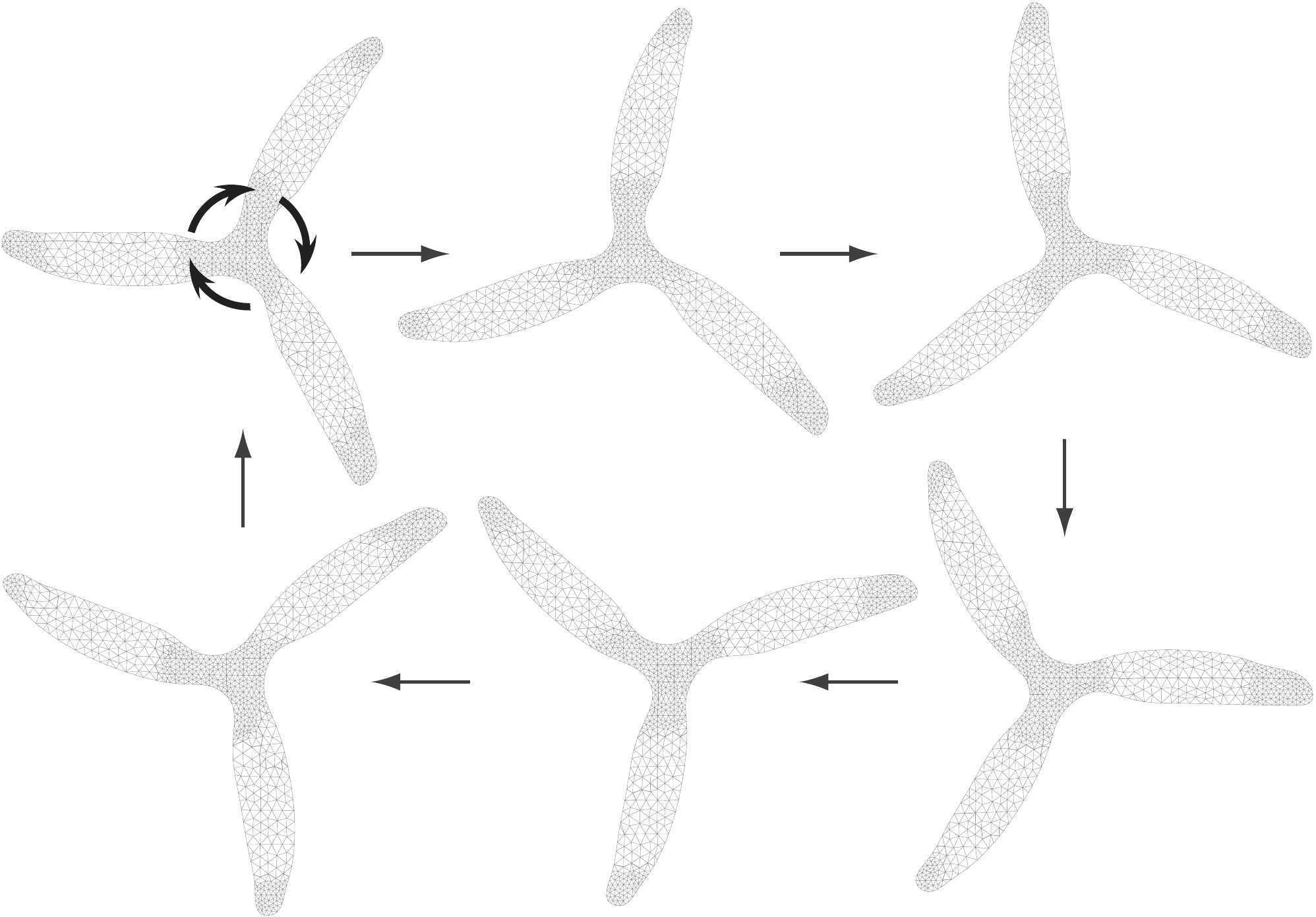}}
  \caption{Example to illustrate the notion of a universal mesh. A
    three blade propeller rotates about an axis perpendicular to its
    plane and passing through its center. It is immersed in a refined
    background mesh of acute angled triangles shown in (a). Using this
    background mesh in the meshing algorithm in Table \ref{algo}
    yields a conforming mesh for each orientation of the propeller; a
    few are shown in (b). Such a background mesh is hence termed a
    universal mesh for the domain of the propeller.}
  \label{fig:propeller}
\end{figure}
Given a triangulation ${\cal T}_h$, let ${\cal D}({\cal T}_h)$ denote
the class of all domains that can be meshed with the algorithm in
Table \ref{algo} using ${\cal T}_h$ as a background mesh. We say that
${\cal T}_h$ is a \emph{universal mesh} for domains in ${\cal D}({\cal
  T}_h)$. The utility of this concept lies in the fact that if
$\{\Omega_t\}_t$ is the time evolution of a domain $\Omega_0\in {\cal
  D}({\cal T}_h)$, then often $\{\Omega_t:0\leq t\leq T\}\subset {\cal
  D}({\cal T}_h)$ for a reasonably large time $T>0$. As the domain
develops small features or undergoes topological changes, it may no
longer belong to ${\cal D}({\cal T}_h)$.

We illustrate this idea with the example shown in
Fig. \ref{fig:propeller}. The domain is a three-blade propeller that
rotates about an axis perpendicular to its plane and passing through
the center. The background mesh shown in Fig. \ref{fig:bg-splines}
consists of only acute angled triangles. It is refined at the center
and along the tips of the blades to resolve the larger curvatures
there.  This mesh yields a conforming mesh for every orientation of
the propeller; a few are shown in
Fig. \ref{fig:propeller-meshes}. This background mesh is hence a
universal mesh for each configuration of the propeller.

An important question in practice is knowing when $\Omega$ belongs to
${\cal D}({\cal T}_h)$. Precisely characterizing ${\cal D}({\cal
  T}_h)$ is presumably very difficult. Fortunately, it is not
essential. Rather, the key step in checking if $\Omega\in {\cal
  D}({\cal T}_h)$ is knowing if ${\cal T}_h$ is sufficiently refined
in the vicinity of $\partial\Omega$. For this, we will require good,
computable and local estimates for the mesh size in order for the
meshing algorithm to succeed. As a step in this direction, in
\cite{rangarajan2011analysis} we provided upper bounds for the mesh
size to guarantee a parameterization of $\partial\Omega$ over positive
edges in ${\cal T}_h$. An analysis of the meshing algorithm is
required to derive a similar bound for the required mesh size of the
background mesh.

\subsection{Details for the implementation}
\label{subsec:details}
An implementation of the meshing algorithm is provided in appendix
\ref{sec:algo}. We discuss a few details here.
\begin{enumerate}[(i)]
\item \textbf{Identifying vertices in $\Omega$:} The first step in
  Table \ref{algo} requires identifying which vertices of ${\cal T}_h$
  lie in $\Omega$. This is simplest when $\Omega$ is represented
  implicitly, as $\Omega = \{x\in{\mathbb R}^2:\Psi(x)<0\}$. For then,
  a vertex $v$ in ${\cal T}_h$ belongs to $\Omega$ \emph{iff}
  $\Psi(v)<0$. If such a level set function $\psi$ is not known
  \emph{a priori}, it can be chosen to be the signed distance function
  to $\partial\Omega$, $\phi:{\mathbb R}^2\rightarrow {\mathbb R}$
  defined as
\begin{align}
  \phi (x) = \begin{cases}
    -\text{distance}(x,\partial\Omega) ~ &\text{if}~x\in \Omega, \\
    \text{distance}(x,\partial\Omega) ~ &\text{otherwise}.
  \end{cases}
  \label{eq:phi}
\end{align}

\item \textbf{Closest point projection:} Mapping vertices of positive
  edges onto $\partial\Omega$ requires computing the closest point
  projection to $\partial\Omega$, $\pi:{\mathbb
    R}^2\rightarrow \partial\Omega$ defined as
\begin{align}
  \pi(x) &= \arg\min_{y\in \partial\Omega} \text{distance}(x,y).
  \label{eq:pi}
\end{align}
For $C^2$-regular domains, $\pi$ is related to the signed distance
function $\phi$ by
\begin{align}
  \pi(x) &= x - \phi(x)\nabla\phi(x) \label{eq:pi-phi}
\end{align}
sufficiently close to $\partial\Omega$, see \cite[Theorem
2.2]{rangarajan2011analysis}. Observe that in the meshing algorithm,
$\pi$ needs to be evaluated only over positive edges and that these
edges are by definition within a distance $h$ from
$\partial\Omega$. Hence relation \eqref{eq:pi-phi} can be used to
compute $\pi$ if ${\cal T}_h$ is sufficiently refined. This in turn
requires computing $\phi$ and its derivatives close to the
boundary. We refer to appendix A in
\cite{rangarajan2011parameterization} for a discussion on computing
$\phi$, $\pi$, and their derivatives for parametric and implicit
representations of $\partial\Omega$.

\item \textbf{Relaxing vertices away from $\partial\Omega$:} In step 4
  in Table \ref{algo}, vertices in $\Omega$ that lie close to
  $\partial\Omega$ are perturbed away from the boundary. While such
  perturbations can be realized in numerous ways, we have adopted the
  map in \eqref{eq:pert}. Close to $\partial\Omega$, $\nabla\phi(x)$
  equals the unit outward normal to $\partial\Omega$ at
  $\pi(x)$. Hence $\mathfrak{p}_h(x)$ indeed perturbs vertices
  \emph{away} from the boundary. By selecting $r$ to be ${\cal
    O}(h)$ in the definition of $\mathfrak{p}_h$, only a small number
  of vertices near the boundary are perturbed. Such a scaling is also
  essential in the definition of $\mathfrak{p}_h$ because $\nabla\phi$
  may be defined only in a small neighborhood of $\partial\Omega$.  In
  our examples, we pick $\eta\simeq 0.3$ and $r \simeq
  3h$. Observe that $\mathfrak{p}_h$ does not move vertices that lie
  on $\partial\Omega$. Hence steps 3 and 4 in Table \ref{algo} move
  exclusive sets of vertices. These two steps can therefore be
  performed in either order.

\item \textbf{Tolerances and Round-off:} In identifying which vertices
  lie in $\Omega$ (step 1 in Table \ref{algo}), the effect of
  tolerances and round-off errors is perhaps unavoidable. As a result,
  a vertex in $\Omega$ may be (mis)identified as lying in $\Omega^c$
  and vice versa. The effect of tolerances can in fact be understood
  as introducing small perturbations in the boundary. Incorrectly
  identifying vertices in $\Omega$ will change the collection of
  positively cut triangles, positive edges and hence the resulting
  mesh for $\Omega$. However, the resulting mesh will be valid
  provided conditioning angles remain acute. In particular, if
  triangles in the vicinity of $\partial\Omega$ are acute angled, the
  choice of tolerances and the effect of round-off errors is not
  critical. The resulting mesh may depend on their choice but will be
  valid nonetheless.

\item \textbf{Background meshes:} For a given curved domain $\Omega$,
  assumptions (b)--(d) in \S\ref{subsec:assumptions} impose
  restrictions on the background mesh ${\cal T}_h$.  Since the polygon
  triangulated by ${\cal T}_h$ is quite arbitrary, assumption (b) is
  easily satisfied. A simple way to satisfy the acute conditioning
  angle requirement (c) is to ensure that triangles in the vicinity of
  $\partial\Omega$ are acute angled. An even simpler way is to use a
  background mesh of all acute angled triangles. For instance, use a
  mesh of all equilateral triangles as done in the example in
  Fig. \ref{fig:bunny}.

  In practice, it is desirable to use an adaptively refined background
  mesh, depending on the geometric features of the boundary or on the
  solution being approximated. A convenient way of doing so is by
  triangulating adaptively refined quadtrees. Bern et
  al. \cite{bern1994provably} provide stencils of acute angled
  triangles to tile quadtrees. The interior angles of triangles in
  these stencils lie between $36^\circ$ and $80^\circ$. Therefore, the
  resulting background meshes automatically satisfy the acute
  conditioning angle requirement. The background mesh shown in
  Fig. \ref{fig:bg-splines} was constructed in this way. We refer to
  \cite{rangarajan2011parameterization} for more examples of
  background meshes constructed from adaptively refined quadtrees.
\end{enumerate}


\section{Exactly conforming curved elements}
\label{sec:curved-elements}
Curvilinear discretizations provide high-order accurate approximations
for curved domains, compared to polygonal ones that results from the
meshing algorithm discussed above. Curved finite elements constructed
using such discretizations are indispensable for optimal accuracy with
high-order interpolations.  Curvilinear discretizations broadly fall
into two categories. In the first kind, curved triangles conform
\emph{exactly} to the domain. In the other, curved triangles
approximate the domain sufficiently well and are usually defined via
isoparametric mappings. We consider the former here and the latter in
\S\ref{sec:isoparametric}.

Constructing mappings from straight to curved triangles, even with a
conforming mesh, is a delicate task because there are two conflicting
requirements. The resulting curved triangle should approximate the
domain well. Yet, it should be a sufficiently small perturbation of
the straight one if interpolation estimates on the latter are expected
to translate into optimal ones over the curved triangle, see
\cite{ciarlet1972interpolation}.  Below we give one such mapping,
which generalizes the ones in \cite{gordon1973transfinite,
  zlamal1973curved} to the case of nonconforming background meshes.

\subsection{Exactly conforming triangles}
\begin{figure}
  \centering
  \subfloat[\label{fig:T1T2}]
  {\includegraphics[scale=0.76]{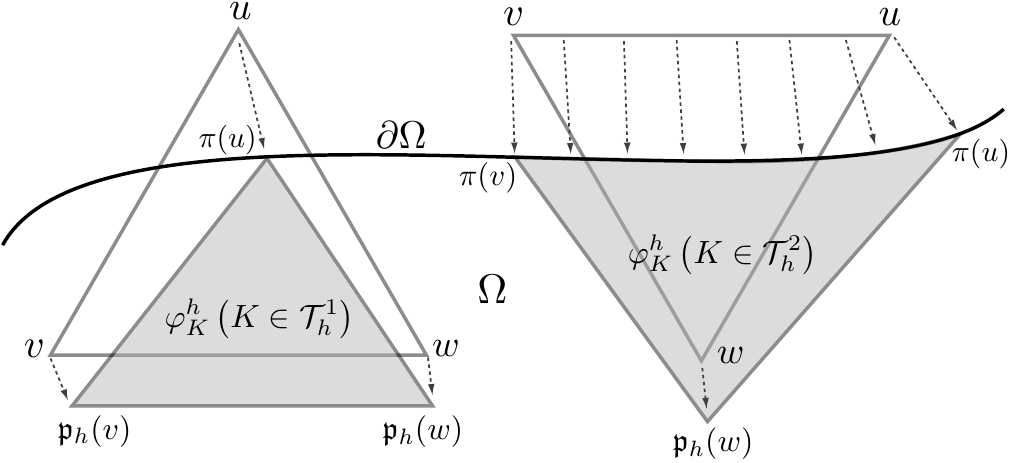}}
  \hfill
  \subfloat[\label{fig:bunny-curved}]
  {\includegraphics[scale=0.45]{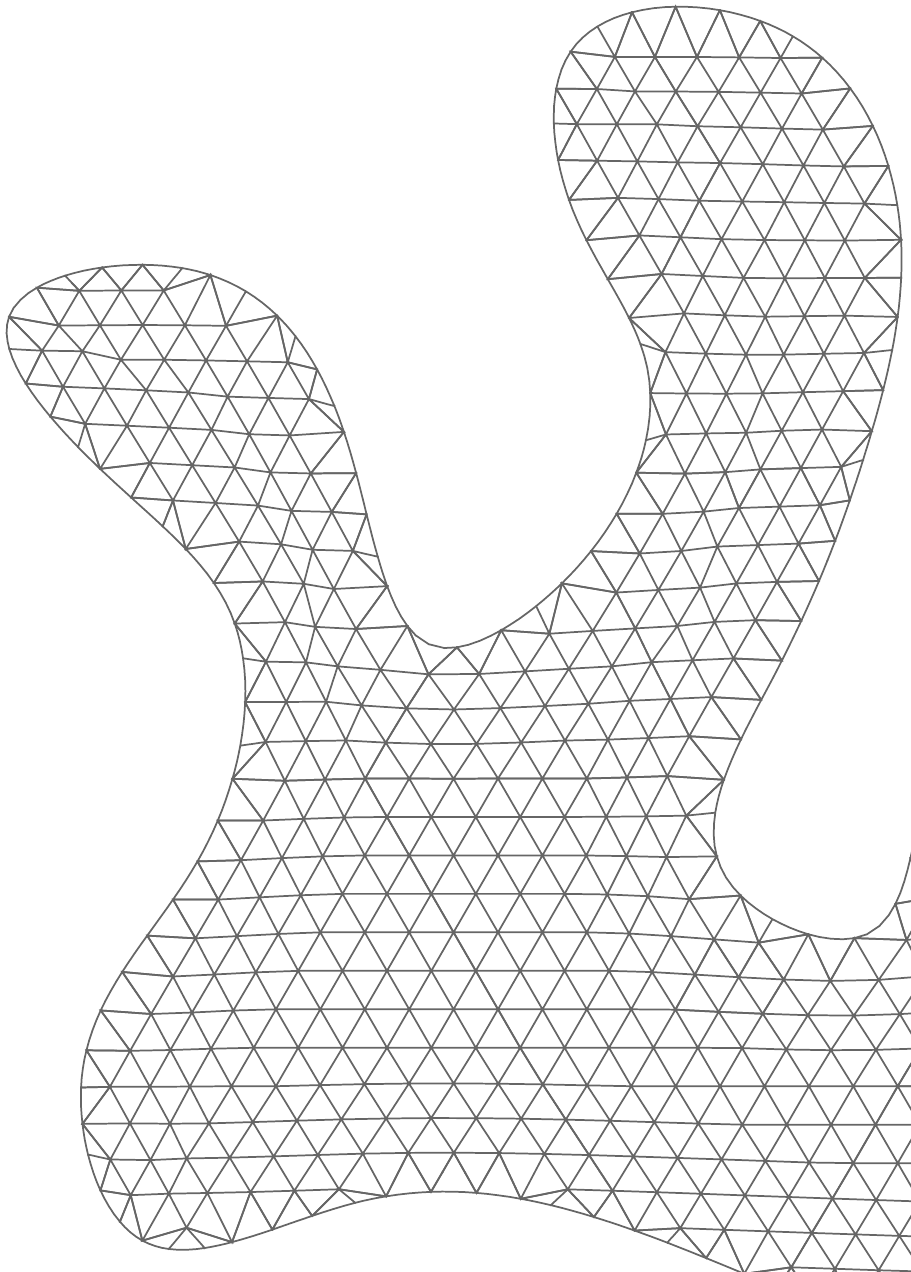}}
  \caption{Mappings from triangles in a background mesh to curved
    triangles that exactly conform to an immersed domain. As shown in
    (a), the mapping $\varphi_K^h$ takes positively cut triangles to
    curved ones that conform to the boundary exactly. This is achieved
    by ensuring that the restriction of $\varphi_K^h$ to a positive
    edge equals the closest point projection $\pi$. Triangles with two
    or more vertices in $\Omega$, i.e., in ${\cal T}_h^{0,1}$, are
    mapped affinely. Figure (b) shows part of the curvilinear mesh
    obtained by using $\varphi_K^h$ in the example in
    Fig. \ref{fig:bunny-a}.}
  \label{fig:curved-mesh}
\end{figure}
Defining a curvilinear mesh that conforms exactly to $\Omega$ requires
only a subtle modification of the meshing algorithm--- instead of
mapping vertices of positive edges onto $\partial\Omega$, we map
positive edges themselves onto $\partial\Omega$. Analogous to the
mapping $M_K^h$ in \eqref{eq:Mk} that defined the meshing algorithm,
we construct a mapping $\varphi_K^h$ triangle-wise over the collection
${\cal T}_h^{0,1,2}$. To this end, consider $K\in {\cal T}_h^{0,1,2}$
with vertices $\{u,v,w\}$ ordered such that $\phi(u)\geq \phi(v)\geq
\phi(w)$. For $K\in {\cal T}_h^{0,1}$, set $\varphi_K^h := M_K^h$.
Over positively cut triangles $K\in {\cal T}_h^2$, define
\begin{align}
  \varphi_K^h(x) &:=
  \frac{1}{2(1-\lambda_u)}\left[\lambda_v\pi\left(\lambda_u\,u +
      (1-\lambda_u)\,v \right) +\lambda_u\lambda_w\pi(u)\right]
  \notag \\
  &\qquad + \frac{1}{2(1-\lambda_v)}\left[\lambda_u\pi\left(
      (1-\lambda_v)\,u + \lambda_v\,v \right)
    +\lambda_v\lambda_w\pi(v)\right]
  \notag \\
  & \quad \qquad + \lambda_w\mathfrak{p}_h(w).
  \label{def:map}
\end{align}
Note that as in \eqref{eq:Mk}, the dependence on $x$ in
\eqref{def:map} is implicit in the barycentric coordinates
$\lambda_u,\lambda_v$ and $\lambda_w$. Unlike $M_K^h$ in \eqref{eq:Mk}
however, $\varphi_K^h$ in \eqref{def:map} is no longer affine over
positively cut triangles. Fig. \ref{fig:T1T2} depicts the action of
$\varphi_K^h$ on triangles in ${\cal T}_h^{1,2}$.
Fig. \ref{fig:bunny-curved} shows part of the curvilinear mesh
obtained by using the map $\varphi_K^h$ in the example in
Fig. \ref{fig:bunny}. Since $\varphi_K^h$ differs from $M_K^h$ only
over triangles in ${\cal T}_h^2$, the curvilinear mesh in
Fig. \ref{fig:bunny-curved} differs from the mesh in
Fig. \ref{fig:bunnymesh-2} only over the $234$ positively cut
triangles.

Let us examine the definition of $\varphi_K^h$ for $K\in {\cal T}_h^2$
in \eqref{def:map}.  By the assumed ordering of vertices, the edge
$\overline{uv}$ joining vertices $u$ and $v$ is the positive edge of
$K$. On this edge $\lambda_w=0$ and $\lambda_u+\lambda_v=1$. So
\begin{align}
  \varphi_K^h(x\in\overline{uv}) 
  &=  \frac{1}{2}\pi((1-\lambda_v)\,u + \lambda_v\,v) 
  + \frac{1}{2}\pi(\lambda_u\,u + (1-\lambda_u)\,v) = \pi(x). \label{eq:phi-1}
\end{align}
Hence $\varphi_K^h$ equals the closest point projection over the
positive edge $\overline{uv}$. This shows that $\varphi_K^h$ maps the
positive edge onto $\partial\Omega$, as depicted in
Fig. \ref{fig:T1T2}. On the edge $\overline{uw}$, $\lambda_v=0$ and
$\lambda_u+\lambda_w=1$. Then \eqref{def:map} reduces to
\begin{align}
  \varphi_K^h(x\in{\overline{uw}}) &=
  \frac{\lambda_u\lambda_w}{2(1-\lambda_u)}\pi(u) + \frac{1}{2}\pi(u)
  + \lambda_w\mathfrak{p}_h(w) = \lambda_u\,\pi(u) +
  \lambda_w\mathfrak{p}_h(w),
  \label{eq:phi-2}
\end{align}
which is an affine map. Similarly,
\begin{align}
  \varphi_K^h(x\in{\overline{vw}}) 
  &= \lambda_v\,\pi(v) + \lambda_w\mathfrak{p}_h(w).
  \label{eq:phi-3}
\end{align}
Eqs.\eqref{eq:phi-1}, \eqref{eq:phi-2} and \eqref{eq:phi-3} show that
$\varphi_K^h$ can be interpreted as the interpolation to
$\overline{K}$, of a map that equals $\pi$ on the positive edge and is
affine on the remaining two. This point of view is also adopted in
\cite{gordon1973transfinite,mansfield1978approximation}. For this
reason, mappings such as $\varphi_K^h$ in \eqref{def:map} are also
commonly termed blending maps and transfinite interpolations.

\textbf{Remark}: We ought to mention an alternate construction for
mapping positively cut triangles to curved ones that explicitly uses
the meshing algorithm as an intermediate step.  In such a
construction, the domain $\Omega$ is first meshed using the algorithm
in \S\ref{sec:meshing} and the resulting mesh is then transformed to a
curvilinear one that conforms to $\Omega$. More precisely, with
$K_S:=M_K^h(K)$, the mapping $\psi_K^h$ defined as
\begin{align}
  \psi_K^h := \begin{cases}
    M_K^h ~&\text{if}~K\in {\cal T}_h^{0,1}, \\
    \varphi_{K_S}^h\circ M_K^h &\text{if}~K\in {\cal T}_h^2 
  \end{cases}
  \label{def:psi}
\end{align}
maps triangles in the collection ${\cal T}_h^{0,1,2}$ to a curvilinear
mesh that conforms exactly to $\Omega$. Of course, the maps
$\varphi_K^h$ and $\psi_K^h$ differ only for positively cut
triangles. For $K\in {\cal T}_h^2$, observe
\begin{enumerate}[(i)]
\item that $\psi_K^h$ first maps $K$ to a conforming triangle $K_S$
  and then transforms $K_S$ to a curved triangle, and
\item that even though
  $\varphi_K^h(\overline{K})=\psi_K^h(\overline{K})$ (as sets in
  ${\mathbb R}^2$), $\varphi_K^h\neq \psi_K^h$ in general. The two
  will however be close in a pointwise sense.
\end{enumerate}
The distinction between the $\varphi_K^h$ and $\psi_K^h$ for
positively cut triangles is illustrated in Fig. \ref{fig:curved-fe}.

The conditions in \S\ref{subsec:assumptions} suffice for the mappings
$\varphi_K^h$ and $\psi_K^h$ to be well defined and bijective.  The
mesh size required near the boundary may be smaller for curvilinear
discretizations of $\Omega$ using $\varphi_K^h$ or $\psi_K^h$ compared
to the mesh size required for the meshing algorithm. This is because
we require the Jacobian of these maps to be positive over the entire
element.

\subsection{High-order finite elements}
\label{subsec:cvg}
\begin{figure}
  \centering
  \includegraphics[scale=1]{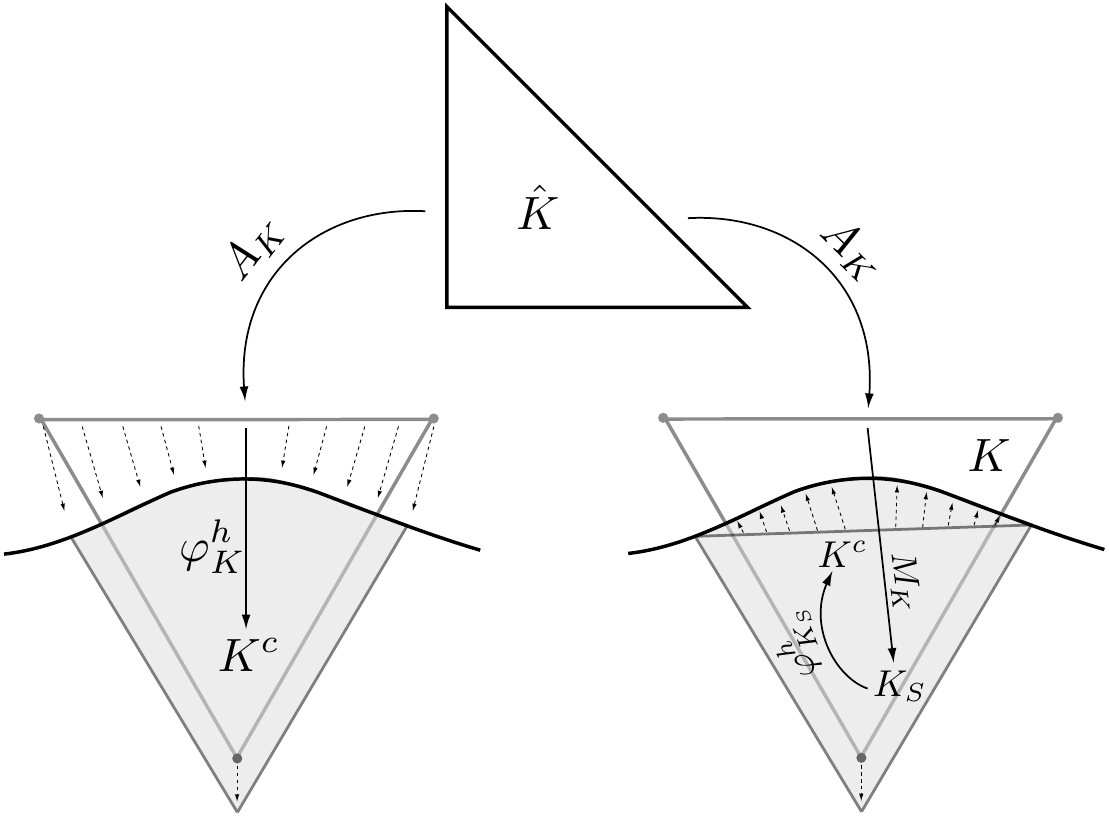}
  \caption{Defining high-order finite elements by constructing
    mappings from the reference element $\hat{K}$ to the curved
    element $K^c$. The figure shows two such mappings. The mapping on
    the left, namely $\varphi_K^h\circ A_K$, first transforms
    $\hat{K}$ to the positively cut triangle $K$ in the background
    mesh and then uses $\varphi_K^h$ to map $K$ to the curved triangle
    $K^c$.  The second construction on the right, given by
    $\varphi_{K_S}^h\circ M_K\circ A_K$ (also equal to $\psi_K^h\circ
    A_K$), uses the meshing algorithm as an intermediate step. The
    reference triangle is mapped to $K$ with $A_K$, then $K$ mapped to
    the triangle $K_S=M_K(K)$ which is finally transformed to $K^c$ with
    the map $\varphi_{K_S}^h$.}
  \label{fig:curved-fe}
\end{figure}
The mappings $\varphi_K^h$ and $\psi_K^h$ from straight to curvilinear
triangles facilitate a natural construction of curved Lagrange finite
elements. The idea behind the construction is illustrated in
Fig. \ref{fig:curved-fe}.

\subsubsection{Curved finite elements}
Introduce the reference triangle $\hat{K}\subset {\mathbb R}^2$ and
the finite element triplet $(\hat{K},\hat{N}^{k},\hat{{\mathbb P}}^k)$. As
usual, $\hat{{\mathbb P}}^k$ is the set of polynomials over $\hat{K}$
of degree at most $k$ and $\hat{N}^{k}=\{\hat{N}_a\}_a$ is the set of
shape functions that constitute a basis for $\hat{{\mathbb
    P}}^k$. Associated with $\hat{N}^{k}$ are the nodes
$\{\hat{z}_a\}_a\in \hat{K}$ which are such that
$\hat{N}_a(\hat{z}_b)=\delta_{ab}$.

Let $A_K:\hat{K}\rightarrow K$ be an affine map from $\hat{K}$ to
$K$. The curved finite element corresponding to $K\in {\cal
  T}_h^{0,1,2}$ derived from $(\hat{K},\hat{N}^{k},\hat{{\mathbb P}}^k)$
is denoted by $(K^c,N^{{k}},P^k)$, where $K^c = \varphi^h_K({K}) =
\varphi_K^h(A_K(\hat{K}))$ and
\begin{align}
  P^k &= \{\hat{p}\circ A_K^{-1}\circ(\varphi_K^h)^{-1}:\:
  \hat{p}\in \hat{{\mathbb P}}^k\}. \label{eq:func-space}
\end{align}
In particular, shape functions $\{N_a\}_a$ over $K^c$ are defined by
the relation $N_a\circ \varphi_K^h\circ A_K =\hat{N}_a$. Nodes
$\{z_a^c\}_a$ in the curved element are located at
$z_a^c=\varphi_K^h(A_K(\hat{z}_a))$. Hereafter, the set
$\{\hat{z}_a\}$ will be chosen so that finite element functions over
$\Omega$ are in $C^0$.

The choice of $\varphi_K^h$ over $\psi_K^h$ to define the curved
element $(K^c,N^{k},P^k)$ above was arbitrary; replacing each instance
of $\varphi_K^h$ by $\psi_K^h$ yields a curved element as well.  In
fact, it may be more convenient to incorporate $\psi_K^h$ into
existing finite element codes based on conforming meshes. In the
following example as well as the ones in \S\ref{subsec:applications},
we have adopted the curved elements based on the map $\varphi_K^h$.

\subsubsection{Optimal convergence: numerical example}
We demonstrate optimal convergence using the curved finite elements
described above with a numerical example. Although the given
construction for curved elements is a standard one, the example helps
show that the mapping $\varphi_K^h$ in \eqref{def:map} satisfies the
conditions in \cite{ciarlet1972interpolation} for optimal
interpolation estimates over the curved element. We consider the model
problem
\begin{subequations}
  \label{eq:model-problem}
  \begin{align}
    \Delta u &= 0~  \text{in}~\Omega=\{r=\sqrt{x^2+y^2}<1\},  \\
    u &= e^y\sin{x}~ \text{on}~\partial\Omega. \label{eq:model-problem-bc}
\end{align}
\end{subequations}
The solution to \eqref{eq:model-problem} is the smooth function
$u(x,y) = e^y\sin{x}$.  The weak form of \eqref{eq:model-problem}
  is to find $u\in H^1_{\partial} = \{v\in
  H^1(\Omega):v\big|_{\partial\Omega}=e^y\sin x\}$ such that
\begin{align}
  \int_{\Omega} \nabla u\cdot \nabla v \,d\Omega &=0 \quad \forall
  v\in H^1_0(\Omega). \label{eq:model-weak}
\end{align}
\begin{SCfigure}
  \centering
  \scalebox{0.65}{\input{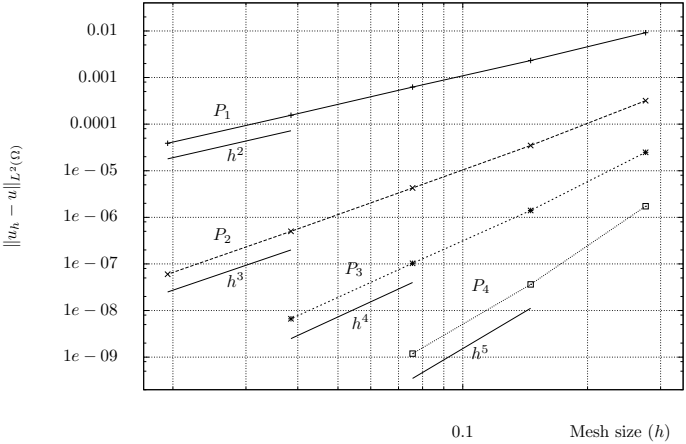}}
   \caption{Optimal convergence of the finite element solution $u_h$,
      computed using exactly conforming curved elements, to the exact
      one $u$ of problem \eqref{eq:model-problem}. The plot shows
      convergence in the $L^2(\Omega)$-norm as the background mesh is
      refined. The rate of convergence is optimal for linear,
      quadratic, cubic and quartic elements. }
  \label{fig:h-cvg}
\end{SCfigure}
To compute finite element approximations $u_h$ of $u$, $\Omega$ is
immersed in background meshes of equilateral triangles. The coarsest
background mesh has mesh size $h_0\simeq 0.27$. Fig. \ref{fig:h-cvg}
shows the convergence of the solution computed with standard Lagrange
elements (over $\hat K$), as the background mesh is refined
($h$-refinement). Dirichlet boundary conditions were imposed by
interpolating the prescribed function in \eqref{eq:model-problem-bc}
at the nodes of curved elements lying on the boundary. We used
sufficiently accurate quadrature rules to evaluate the stiffness
matrix, see \S\ref{subsec:quadrature}. The convergence rate in the
$L^2(\Omega)$-norm is optimal for linear, quadratic, cubic and quartic
elements ($k=1,2,3$ and $4$ respectively).

Examining $\|u-u_h\|_{L^2(\Omega)}$ in Fig. \ref{fig:h-cvg} also
reveals that for a given background mesh, the error decreases with the
element order $k$. This demonstrates that the curved elements are well
suited for $p$-refinement--- progressively accurate solutions can be
computed by just increasing the element order while using the
\emph{same background mesh}.  Moreover, the data in
Fig. \ref{fig:h-cvg} shows that the error is ${\cal O}(h^{k+1})$,
which is optimal in the element order $k$ for each given (sufficiently
small) mesh size $h$ of the background mesh. Such optimal convergence
rates would also be obtained with the isoparametric curved elements
described subsequently in \S\ref{sec:isoparametric}.

\subsection{Applications: Evolving fluid domains}
\label{subsec:applications}
Next, we present two applications using the curved elements described
above.  In both examples, a fixed background mesh serves as the
universal mesh for an evolving fluid domain.

\subsubsection{Flow with a rotating component}
\label{subsec:propeller}
\begin{figure}
\centering
\subfloat[\label{fig:prop-1-1}]{\includegraphics[scale=0.5]{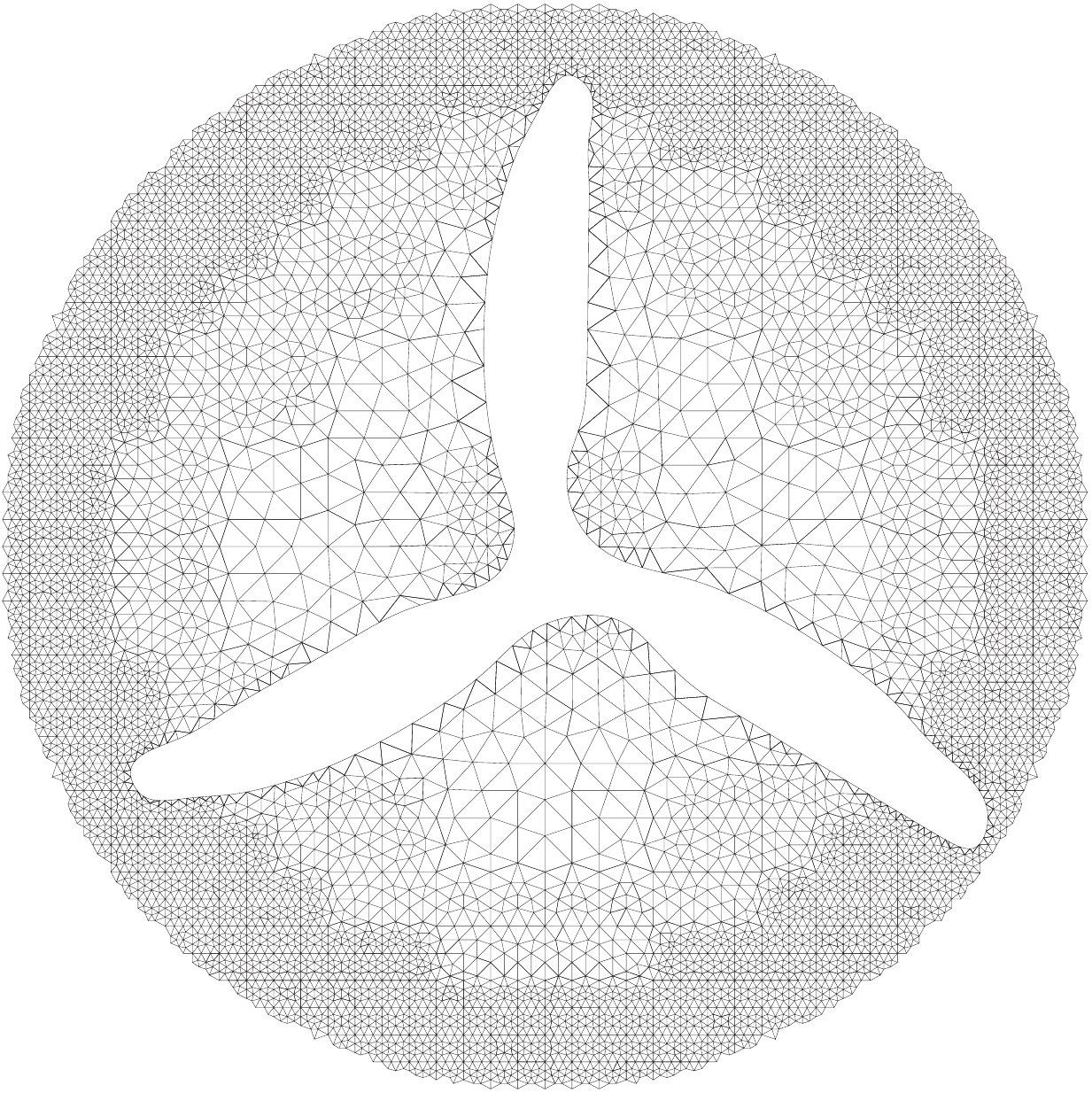}}
\hfill
\subfloat[\label{fig:prop-3-1}]{\includegraphics[scale=0.5]{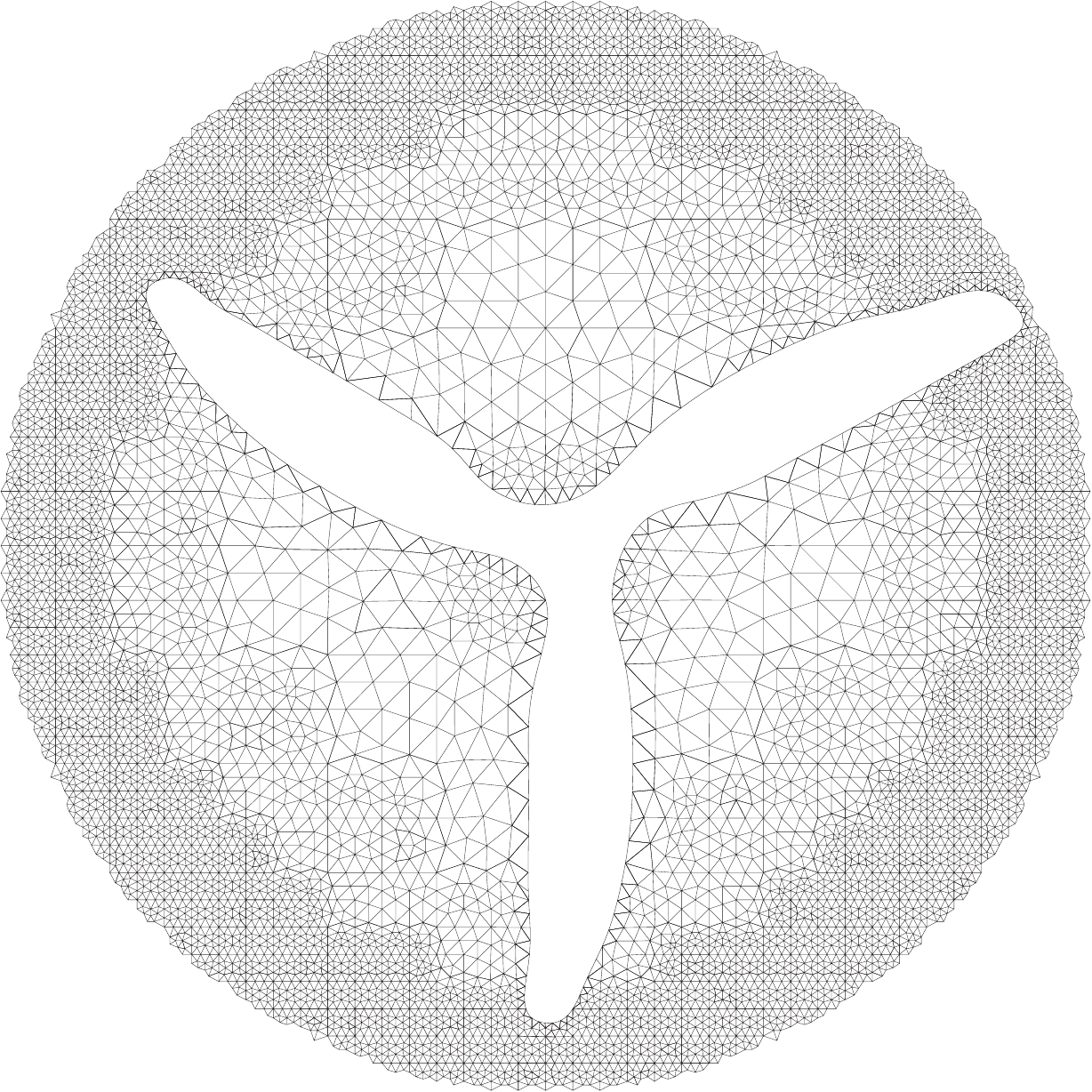}}
\\
\subfloat[\label{fig:prop-1-2}]{\includegraphics[scale=0.6]{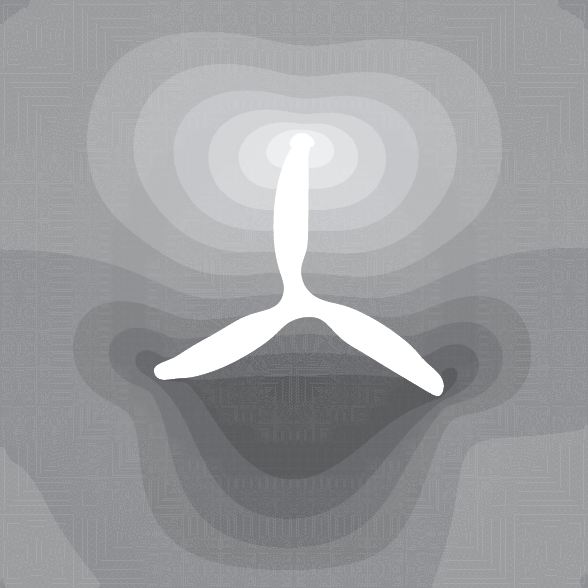}}
\hfill
\subfloat[\label{fig:prop-3-2}]{\includegraphics[scale=0.6]{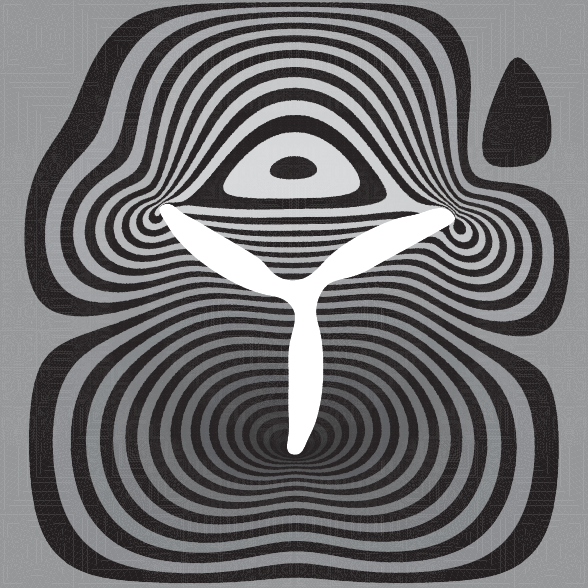}}
\caption{Simulating Stokes flow driven by a rotating propeller in a
  closed container. At each instant, the same background mesh is used
  to determine a curvilinear mesh that conforms exactly to the fluid
  domain. Figures (a) and (b) show the elements close to the propeller
  in the resulting mesh at two distinct times. Contours of the
  horizontal component of the flow velocity at these times are shown
  in figures (c) and (d). To highlight the flow pattern, the contours
  are shown with zebra shading in (d). }
\end{figure}
We consider the example mentioned in \S\ref{sec:introduction}, of a
propeller mixing fluid in a closed container. The problem setup is the
same one illustrated in Fig. \ref{fig:propeller}, although the
container $B$ is larger. The propeller $P$ is assumed to be rigid and
impermeable. It rotates with constant angular velocity $\omega$ about
an axis passing through its center and perpendicular to its plane.
The fluid in the container is incompressible and has viscosity
$\mu$. Its kinematics is governed by the familiar equations for Stokes
flow,
\begin{subequations}
  \label{eq:stokes-1}
\begin{align}
  \mu\,\text{div}\left(\nabla{\bf u}_t\right)
  &= \nabla p_t, \label{eq:stokes-1a} \\
  \text{div}\left({\bf u}_t\right) &=0, \label{eq:stokes-1b}
\end{align}
\end{subequations}
relating the flow velocity ${\bf u}_t$ and pressure $p_t$ at time
$t$. No-slip boundary conditions along the walls of the container and
the boundary of the propeller imply 
\begin{align}
  {\bf u}_t = 
  \begin{cases}
    0 ~ ~ &\text{on}~\partial B, \\
    r\omega \,{\bf e}_{\theta}&\text{on}~\partial P_t,
  \end{cases}
   \label{eq:stokes-2}
 \end{align}
 where $P_t$ is the configuration of the propeller at time $t$ and
 $\{{\bf e}_r,{\bf e}_\theta\}$ is a canonical polar basis for a polar
 coordinate system with
 origin at the center of $P$. Since \eqref{eq:stokes-1} and
 \eqref{eq:stokes-2} determine the pressure only up to a constant,
 $p_t$ is assigned to be zero at one point in the flow domain, i.e.,
 we set $p_t=0$ at $x_0\in B\setminus P_t$.

 A triangulation of $B$, similar to the one shown in
 Fig. \ref{fig:bg-splines}, serves as the universal mesh for the fluid
 domain $B\setminus P_t$.  The mesh is refined further near the tips
 of the blades to resolve features of the flow there. We adopt
 (curved) Taylor-Hood elements for the finite element solution of
 \eqref{eq:stokes-1}, i.e., the element $(K^c,N^{2},P^2)$ for the
 velocity ${\bf u}_t$ and $(K^c,N^{1},P^1)$ for the pressure $p_t$.
 See \cite{girault1986finite} for a discussion on the Taylor-Hood
 element and for the weak form of this problem, which we have omitted
 here. Dirichlet boundary conditions are imposed by interpolating
 \eqref{eq:stokes-2} at the nodes lying on the boundary.  Figs.
 \ref{fig:prop-1-1} and \ref{fig:prop-3-1} show the curvilinear mesh
 conforming to the fluid domain at three different time
 instants. Since the mesh is quite refined, only the elements near the
 propeller are shown. Corresponding to these orientations of the
 propeller, Figs. \ref{fig:prop-1-2} and \ref{fig:prop-3-2} show
 contours of the horizontal component of the velocity computed with
 $\mu=0.01$ and $\omega=2$.

 \subsubsection{Flow interaction with a rigid disc}
\label{subsec:stokes}
\begin{figure}
  \centering
  \subfloat[\label{fig:ball-initial}]{\includegraphics[scale=0.29]{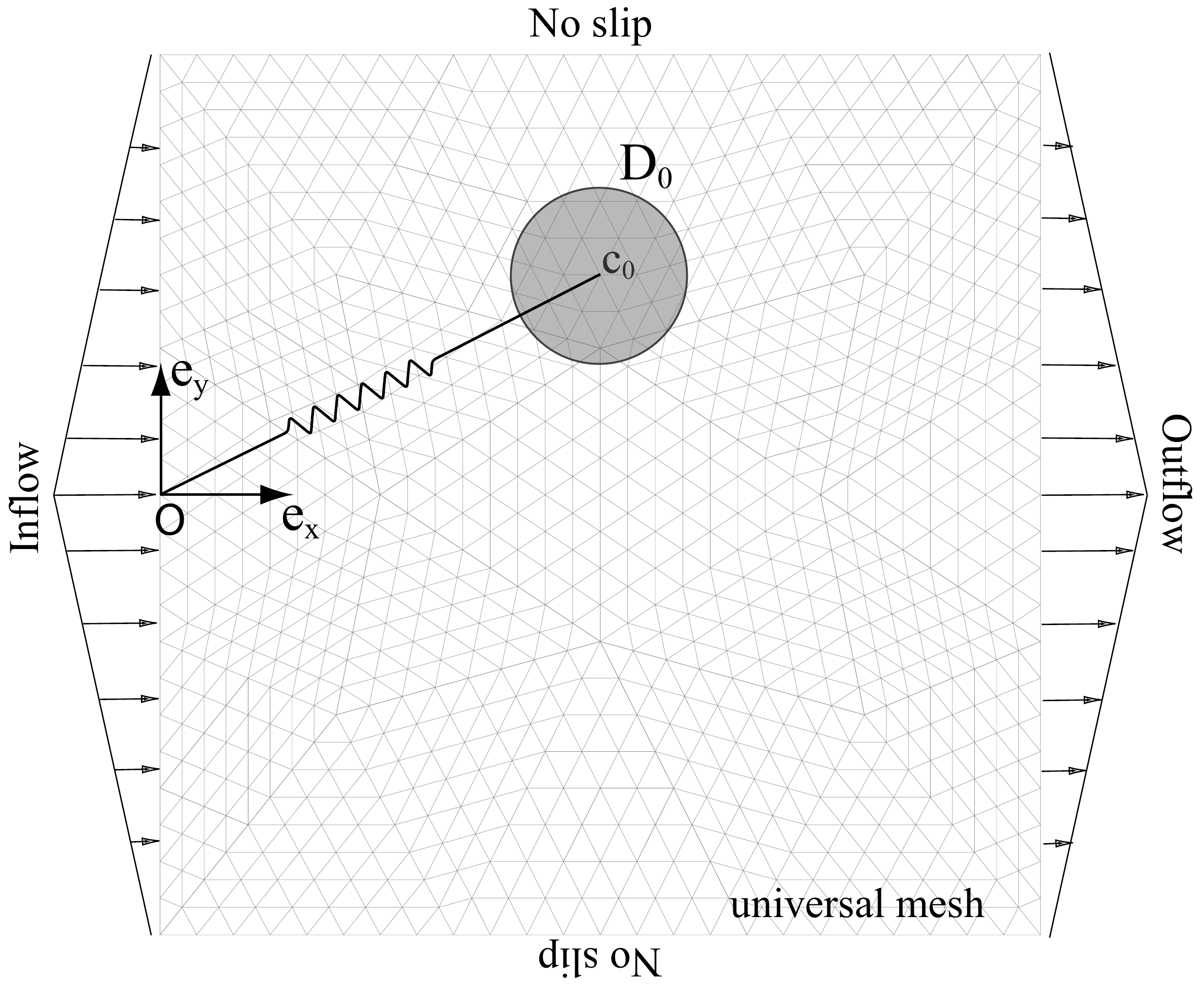}}
  \hfill
  \subfloat[\label{fig:ball-final}]{\includegraphics[scale=0.31]{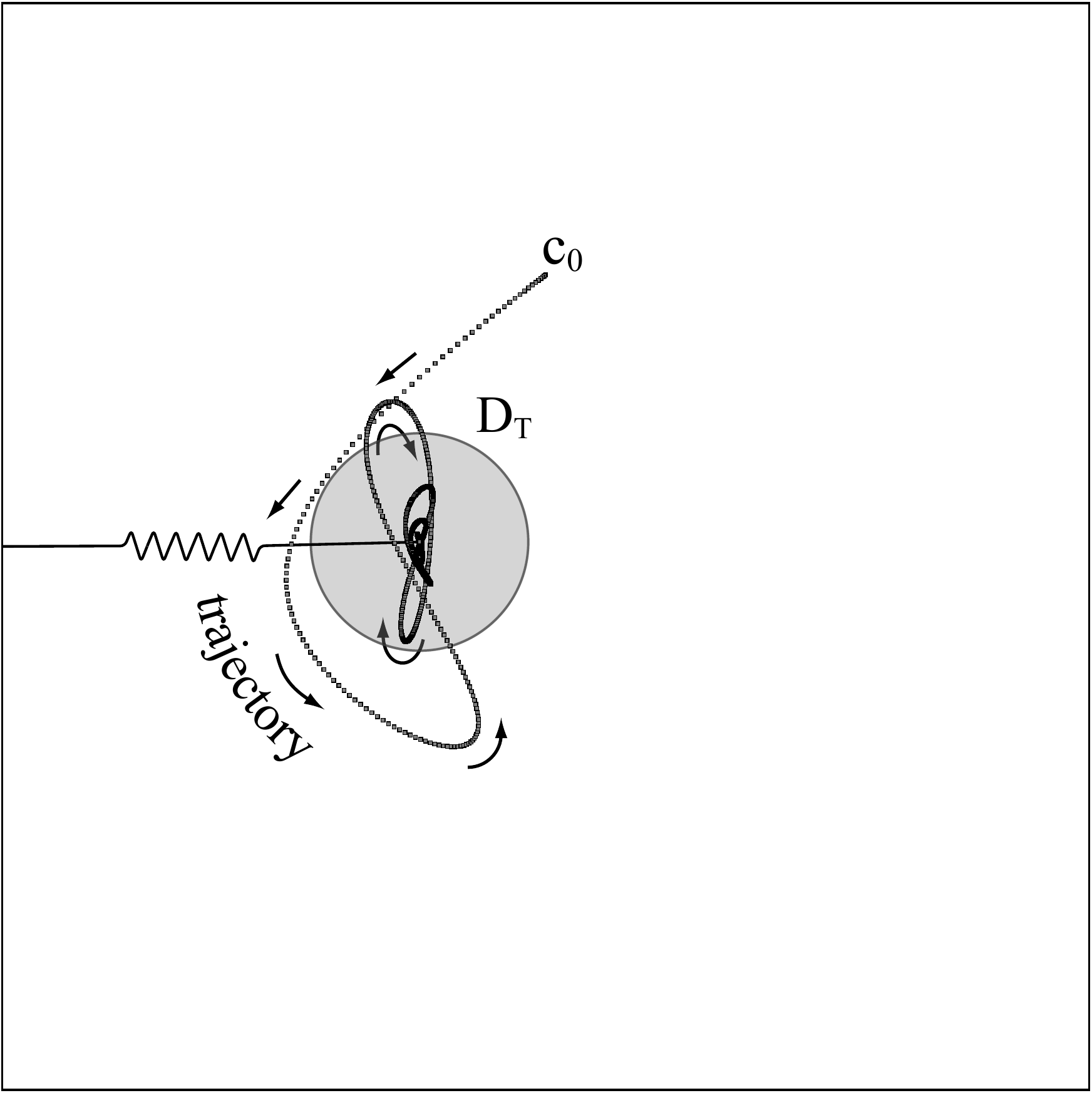}}
  \caption{Figure (a) shows the initial setup for the problem of a
    rigid disc $D$ interacting with an incompressible fluid. The disc
    is attached to the origin by a linear spring. Inflow and outflow
    boundary conditions for the flow through the channel $B$ are
    indicated. A simple, unstructured mesh of $B$ serves as the
    universal mesh for the fluid for the entire duration of the
    simulation. Figure (b) shows the trajectory computed for the disc
    as it moves from its initial position to an (approximate)
    equilibrium position.}
\label{fig:ball-initial-final}
\end{figure}
In the second example, we consider the interaction between a fluid and
a rigid solid. The problem is to determine the trajectory of a rigid
disc $D$ of radius $R$ and mass $m$ immersed in an incompressible,
viscous fluid flowing through a square shaped channel $B$ of side
$L$. The disc is attached to the origin $O$ located at the mid-point
of the left end of the channel by a linear spring with spring constant
$k$ and equilibrium length $\ell_0$. The problem setup is shown in
Fig. \ref{fig:ball-initial}.

We assume that the kinematics of the fluid is governed by the
equations for Stokes flow given in \eqref{eq:stokes-1}, retaining the
notation introduced there. In a Cartesian coordinate system  centered at $O$, inflow and outflow boundary
conditions are prescribed at the two ends of the flow channel as
\begin{align}
  {\bf u}_t &= (L-2|y|){\bf e}_x ~ ~ \text{if}~x = 0, L. \label{eq:bc-1}
\end{align}
Denote the position of the center of the disc at time $t$ by ${\bf
  c}(t)$ and the disc centered at ${\bf c}(t)$ by $D_t$. No-slip
boundary conditions along the horizontal walls of the channel and
along the boundary of the disc imply
\begin{align}
  {\bf u}_t = \begin{cases}
    0 ~ ~ &\text{if}~|y|=L/2, \\
    \dot{\bf c}(t) &\text{on}~\partial D_t.
  \end{cases}
  \label{eq:bc-2}
\end{align}
Force balance for the disc is given by 
\begin{align}
  m\,\ddot{{\bf c}} &= k\left(1-\frac{\ell_0}{|{\bf c}|}\right){\bf c} 
  + \int_{\partial D_t}{\boldsymbol \sigma}_f\cdot {\bf n}_t\,ds,
  \label{eq:force-solid}
\end{align}
where ${\bf n}_t$ is the unit outward normal to $\partial D_t$ and the
stress ${\boldsymbol \sigma}_f$ in the fluid is computed as
\begin{align}
  {\boldsymbol \sigma}_f &= -p_t\,{\mathbb I} 
  + \mu\left(\nabla{\bf u}_t+\nabla{\bf u}_t^T\right). \notag
\end{align}

\begin{SCfigure}[0.5]
  \centering
  \includegraphics[scale=0.35]{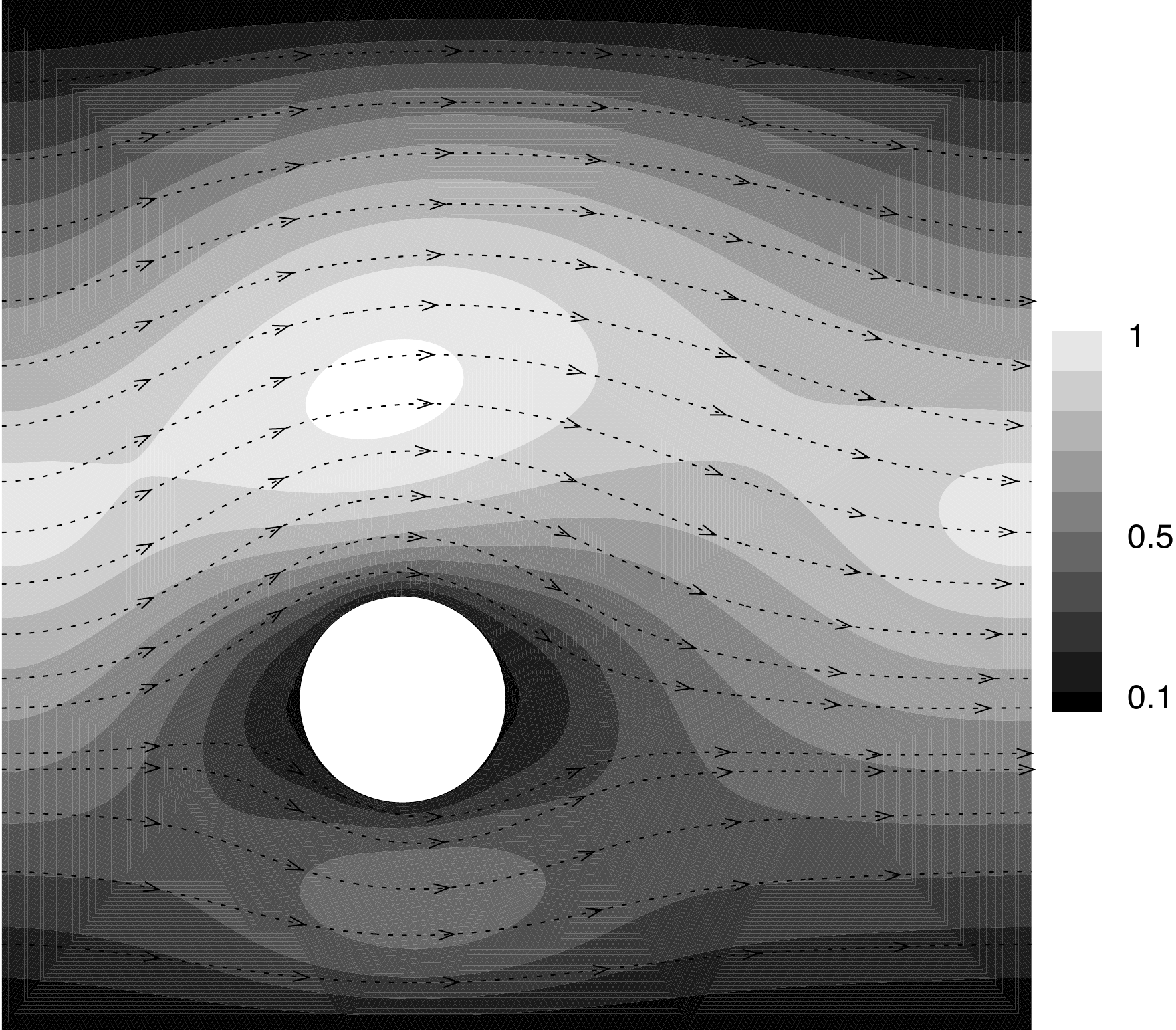}
  \hfill
  \caption{Contours of the horizontal component of the flow velocity
    at a non-equilibrium position of the disc in the problem described
    in Fig. \ref{fig:ball-initial}.}
  \label{fig:ball-Vx}
\end{SCfigure}

Balance equations \eqref{eq:stokes-1}, \eqref{eq:force-solid},
boundary conditions \eqref{eq:bc-1}, \eqref{eq:bc-2}, initial
conditions ${\bf c}(0)={\bf c}_0$, $\dot{{\bf c}}(0)=0$, and
$p_t(x_0\in B\setminus D_t)=0$ together constitute a coupled system of
equations for the unknowns $({\bf u}_t,p_t)$ and ${\bf c}(t)$. We use
(curved) Taylor-Hood elements for the flow solution as before and
adopt a staggered time integration scheme. At each time instant $t_n$,
given the center of the disc ${\bf c}_n$ and its velocity $\dot{{\bf
    c}}_n$, we define curvilinear elements for the flow variables over
$B\setminus D_t$. We then compute the flow solution $({\bf
  u}_n^h,p_n^h)$ at this time. The net force on the disc is evaluated
using \eqref{eq:force-solid}. Using central differences, we update the
position and velocity of the disc to the next instant $t_{n+1}$ and
repeat the process.

The background mesh shown in Fig. \ref{fig:ball-initial} serves as a
universal mesh for the flow domain $B\setminus
D_t$. Fig. \ref{fig:ball-final} shows the trajectory determined for the
disc and its final configuration computed with parameters $L=1$ for
the container, $\mu=0.01$ for the fluid, $R=L/10, m=1$ for the disc,
$k=1, \ell_0=L/4$ for the spring, ${\bf c}_0=(0,L/4)$ for the initial
position of the disc and a time step $\Delta t = 0.05$. At an
intermediate position of the disc, contours of the horizontal
component of the flow velocity along with a few stream lines are shown
in Fig. \ref{fig:ball-Vx}.
The trajectory of the disc plotted in Fig. \ref{fig:ball-final} shows
that the disc eventually settles to an equilibrium position balancing
the forces exerted by the fluid and the spring.

We conclude this section mentioning that there are various numerical
methods in the literature for such problems over changing flow
domains, see for instance \cite{bazilevs2008isogeometric,
  donea2004arbitrary, saksono2007adaptive, wagner2001extended}.

\section{Isoparametric mappings}
\label{sec:isoparametric}
\begin{figure}
  \centering
  \subfloat[\label{fig:P2}]{\includegraphics[scale=0.52]{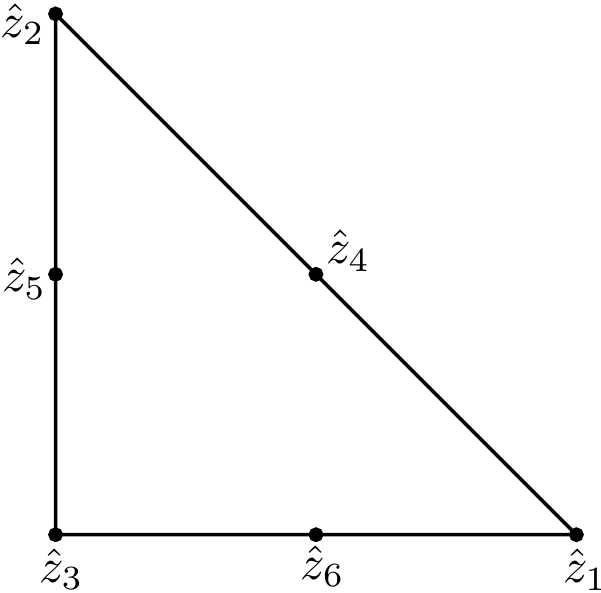}}
  \hfill
  \subfloat[\label{fig:isoparam-1}]{\includegraphics[scale=0.77]{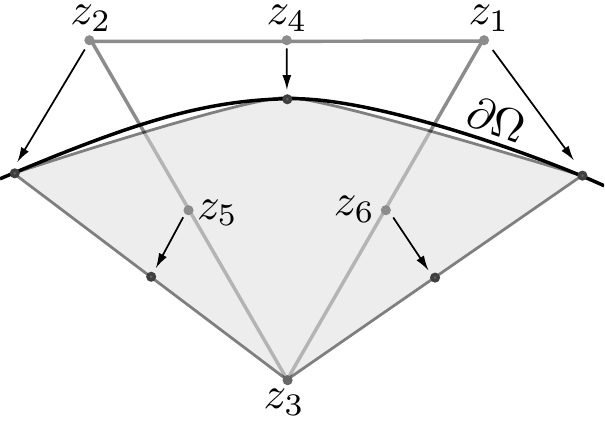}}
  \hfill
  \subfloat[\label{fig:isoparam-2}]{\includegraphics[scale=0.77]{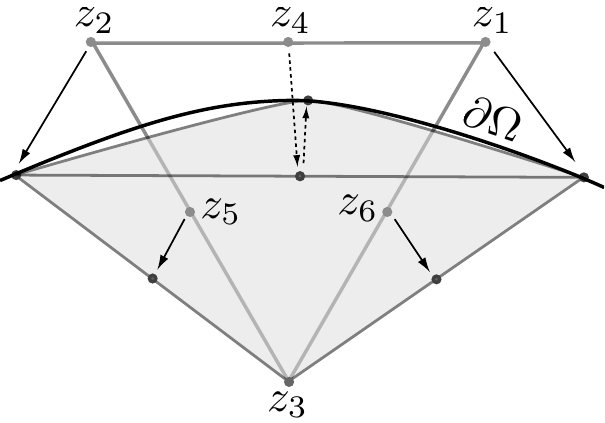}}
  \caption{Isoparametric mappings for positively cut triangles in the
    background mesh. The reference quadratic element is shown in
    (a). The isoparametric maps constructed by interpolating the
    mappings $\varphi_K^h\circ A_K$ and $\psi_K^h\circ A_K$ are shown
    in (a) and (b). The curved boundary of the element interpolates
    the boundary at the points where the nodes $\hat{z}_1,\hat{z}_2$
    and $\hat{z}_4$ are mapped. Notice that the mid-side node
    $\hat{z}_4$ is generally mapped differently by $\varphi_K^h$ and
    $\psi_K^h$. }
  \label{fig:isoparam}
\end{figure}

Isoparametric mappings provide a convenient way of approximating
curved domains with a desired accuracy. A systematic definition of
these polynomial mappings results naturally from interpolating
$\varphi_K^h\circ A_K$ or $\psi_K^h\circ A_K$ at selected
points. Since $\varphi_K^h\circ A_K$ and $\psi_K^h\circ A_K$ are
affine for $K\in {\cal T}_h^{0,1}$, these maps differ from their
interpolants only for positively cut triangles $({\cal T}_h^2)$. With
an isoparametric map, positive edges are mapped to curved ones that
interpolate the boundary at a few points. Fig. \ref{fig:isoparam}
depicts this for the case of a quadratic element.

Following the notation  in \S\ref{subsec:cvg}, introduce the
interpolation operator $\hat{\Pi}^k:f\in [C^0(\hat{K})]^2\rightarrow
\sum_{a}f(\hat{z}_a)\hat{N}_a \in {\mathbb P}^k$. The isoparametric
map over $\hat{K}$ corresponding to triangle $K\in {\cal T}_h^{0,1,2}$
is defined as 
\begin{align}
I_K^h := \hat{\Pi}^k(\varphi_K^h\circ A_K) = \sum_a
\varphi_K^h(z_a)N_a.
\label{eq:isoparam1}
\end{align}
Using $\psi_K^h$ instead of $\varphi_K^h$ yields a different
isoparametric map
\begin{align}
  J_K^h:= \hat{\Pi}^k(\psi_K^h\circ A_K) = \sum_a
  \varphi_{K_S}^h(M_K(z_a))N_a.
  \label{eq:isoparam2}
\end{align}
As illustrated in Fig. \ref{fig:isoparam} for the case of a quadratic
element, $I_K^h\neq J_K^h$ in general. In the figure for instance,
$I_K^h(\hat{z}_4)= \pi((z_1+z_2)/2)$ while $J_K^h(\hat{z}_4) =
\pi((\pi(z_1)+\pi(z_2))/2)$. Nonetheless, the two maps will be close
for small values of $h$. Curved finite elements using these
isoparametric maps are constructed just as in \S\ref{subsec:cvg} by
replacing maps $\varphi_K^h$ or $\psi_K^h$ by their respective
interpolants.

Compared to the exactly conforming curved elements, isoparametric
elements require fewer evaluations of $\pi$ in general. Notice from
\eqref{def:map} that once $\pi$ is computed at the vertices of the
positive edge, defining $I_K^h$ (or $J_K^h$) requires evaluating $\pi$
at most twice \emph{per node} in the element. In contrast, computing
$\varphi_K^h$ (or $\psi_K^h$) at \emph{each quadrature point} requires
two new evaluations in the conforming curved element. Furthermore,
computing derivatives of shape functions in the isoparametric element
does not require computing derivatives of $\pi$. However in the
conforming curved element, shape function derivatives depend on
$\nabla\varphi_K^h$ (or $\nabla\psi_K^h$) which in turn depend on
$\nabla\pi$.

\subsection{Quadrature for curved elements}
\label{subsec:quadrature}
For optimal convergence and accuracy of numerical solutions computed
using curved elements, we naturally require sufficiently accurate
quadrature rules for integration over curvilinear domains. Following
standard practice, these quadrature rules need only be defined over
the reference triangle, since integrals over a curved element $K^c$
can be performed over $\hat{K}$ using the correspondence provided by
the mappings $I_K^h$ (or $J_K^h)$.  We adopted for curved elements the
same quadrature rules needed for straight elements, as explained in
\cite{ciarlet1978finite}. For example, a quadrature rule that exactly
integrates quadratic polynomials over the reference element suffices
for isoparametric quadratic elements, so three quadrature points were
adopted for those. We have used such integration rules in all of our
examples, including the ones with exactly conforming curved
elements. These examples suggest that the quadrature rules for
straight elements are also enough to obtain optimal convergence rates
with exactly curved elements.  We also note that different rules can
be used for integrating each term in the weak form of a problem, for
instance the mass and stiffness matrices, and the force vector.

\subsection{Circular plate in bending}
\label{subsec:plate}
We consider the problem of a thick, circular, elastic plate bending
under the action of a uniform external load. As pointed out in
\cite{babuska1990plate}, this problem highlights the importance of
representing curved boundaries accurately. We present this example to
emphasize the distinction between representing a curved boundary
exactly using the elements in \S\ref{sec:curved-elements} and
approximating it with isoparametric elements as defined
above. Consider a Cartesian coordinate
system $(x,y,z)$ with basis $\{{\bf e}_x,{\bf e}_y,{\bf e}_z\}$. The domain of the problem is the set $\Omega\times
(-t/2,t/2)$, where $\Omega$ is a circle of radius $R=3.142$ centered
at the origin and contained in the plane of ${\bf e}_x,{\bf e}_y$, and
$t$ is the thickness along ${\bf e}_z$. We assume that the
displacement field ${\bf u}$ of the plate is of the form
\begin{align}
  {\bf u}(x,y,z) &= -z\left(\vartheta_x(x,y){\bf
      e}_x+\vartheta_y(x,y){\bf e}_y\right) + w(x,y){\bf
    e}_z, \label{eq:plate}
\end{align}
which corresponds to the Reissner-Mindlin model for a thick plate in
bending, see \cite{babuska1990plate,chapelle2003finite}. In
\eqref{eq:plate}, $w$ is the transverse displacement of points in the
mid-plane $\Omega$ while $\vartheta_x$ and $\vartheta_y$ represent the
infinitesimal rotations of fibers normal to the mid-plane about the
axes ${\bf e}_y$ and ${\bf e}_x$, respectively. We consider a
``soft-simple support'' for the plate, which implies the boundary
conditions
\begin{align}
  \left. 
    \begin{aligned}
      w &= 0 \\
      {\boldsymbol \vartheta}\cdot {\bf t} &= 0
    \end{aligned} \right\} ~\text{on}~\partial\Omega
  \label{eq:plate-bc}
\end{align}
where ${\bf t}$ is the unit tangent to $\partial\Omega$ and $\btheta =
(\vartheta_x,\vartheta_y)$. The plate is loaded by a constant force
$2p$ normal to its top face $\Omega\times \{t/2\}$.  The elasticity
problem is then to find
\begin{align*}
{ \bf u}&\in \left\{-z\,\btheta+w\,{\bf e}_z\,:\,\btheta\in {\bf H}^1_t,w\in
H^1_0(\Omega)\right\}, \\
\text{where}~ ~ 
{\bf H}^1_t(\Omega) &:= \{\btheta\in [H^1(\Omega)]^2\,:\,\btheta\cdot {\bf
  t}=0~\text{on}~\partial\Omega\},
\end{align*} 
which minimizes the strain energy functional
\begin{align}
  I[{\bf u}] &=
  \frac{1}{2}\int_{\Omega\times(-t/2,t/2)}\left\{\lambda [\text{tr}({\boldsymbol
        \varepsilon}({\bf u}))]^2 + 2\mu\,{\boldsymbol
      \varepsilon}({\bf u}):{\boldsymbol \varepsilon}({\bf u})\right\}
  - \int_\Omega p\,w,
  \label{eq:elasticity}
\end{align}
where $\lambda,\mu$ are material parameters called Lam\'{e} constants
and $\varepsilon({\bf u})=(\nabla{\bf u}+{\nabla{\bf u}}^T)/2$ is the
usual infinitesimal strain tensor. Introducing assumption \eqref{eq:plate}
in \eqref{eq:elasticity} and integrating along the thickness reduces
\eqref{eq:elasticity} to a problem over $\Omega$: find $(\btheta,w)\in
{\bf H}^1_t(\Omega)\times H^1_0(\Omega)$ that minimizes the functional
\begin{align}
  F[(\btheta, w)] &=
  \frac{1}{2}\int_{\Omega}\left\{\lambda\left[\text{tr}\left({\boldsymbol \varepsilon}(\btheta) \right)\right]^2 + 2\mu\,{\boldsymbol \varepsilon}({\boldsymbol \vartheta}):{\boldsymbol \varepsilon}({\boldsymbol \vartheta}) \right\} \notag \\
  &\qquad\qquad\qquad\qquad + \frac{6\mu}{t^2}\int_\Omega
  \|{\boldsymbol \vartheta}-\nabla w\|^2 - \frac{12}{t^3}\int_\Omega p
  w. \label{eq:plate-eqn}
\end{align}

We compare the transverse displacements at the center of the plate as
the background mesh for $\Omega$ is refined, while using curved
quadratic elements. We pick $\lambda=\mu=1,t=R/4$ for the plate and
$p=1\e{-3}$ for the loading.  We take as a reference value
$w_0=5.3407075\e{-2}$, computed with exactly conforming curved quartic
elements (see \S\ref{subsec:cvg}) and a refined background mesh of
equilateral triangles.  Table \ref{table:plate} lists the
displacements computed with exactly conforming quadratic elements and
with isoparametric quadratic elements.
\begin{table}
  \caption{Transverse displacements at the center of the circular
    plate $\Omega$ computed 
    with curved quadratic elements. The reference value is
    $w_0=5.3407075\e{-2}$. The columns titled `conforming' and
    `isoparametric'  list the values computed with exactly 
    conforming elements and
    isoparametric elements, respectively. The column `modified isoparametric'
    contains the values computed with isoparametric elements
    but while imposing the constraint $\btheta\cdot {\bf t}=0$ using
    one of the two possible tangents ${\bf t}$ at vertices on the
    curved edges of the element. The coarsest mesh size of the
    background mesh is $h_0\simeq 0.28R$.}
  \begin{center}
    \footnotesize
    \begin{tabular}{|c|c|c|c|}
      \hline
      mesh size & conforming $(\times 10^{-2})$ & isoparametric
      $(\times 10^{-2})$ & modified isoparametric $(\times 10^{-2})$ \\
      \hline
      $h_0$ & $5.3370662$ & \cellcolor[gray]{0.8} $2.2343001$   & \cellcolor[gray]{0.9}$5.3383666$  \\
      $h_0/2$           & $5.3401468$ &  \cellcolor[gray]{0.8} $1.9557577$ & \cellcolor[gray]{0.9}$5.3400195$  \\
      $h_0/4$           & $5.3406606$  &  \cellcolor[gray]{0.8} $1.8019916$ & \cellcolor[gray]{0.9}$5.3406438$  \\
      $h_0/8$           & $5.3407108$ &  \cellcolor[gray]{0.8} $1.6864489$  & \cellcolor[gray]{0.9}$5.3407156$   \\
      \hline
    \end{tabular}
  \end{center}
  \label{table:plate}
\end{table}

From the table, we see that the displacements computed with the
conforming elements converge to $w_0$. But somewhat surprisingly,
those computed with quadratic isoparametric elements fail to even come
close to $w_0$. This is a consequence of enforcing the constraint on
rotations in \eqref{eq:plate-bc} on the approximate curved boundary
realized with isoparametric elements. The unit tangent to this
boundary fails to be continuous at vertices that lie on
it. Consequently, rotations equal zero at each vertex on the
boundary. With exactly conforming curved elements, $\partial\Omega$ is
represented exactly and this issue is avoided.

The above discussion shows that the constraint in \eqref{eq:plate-bc}
needs to be enforced differently. For instance, we could select one of
the tangents at each vertex on the approximate boundary to enforce the
constraint $\btheta\cdot {\bf t}=0$. The displacements computed with quadratic
isoparametric elements by enforcing the constraint in this way are listed under the column title `modified
isoparametric' in Table \ref{table:plate}. These values are clearly
more accurate and converge to $w_0$.
 
\noindent
\textbf{Remark}: We have used the same finite element spaces for both
transverse displacements ($w$) and rotations
($\vartheta_x,\vartheta_y$) in the above calculations. It is well
known that for thin plates ($t\ll R$), such a choice of spaces in the
Reissner-Mindlin model results in locking,
cf. \cite{chapelle2003finite}. To avoid adopting very small mesh sizes
for accuracy, we deliberately chose a large thickness $t=R/4$ in the
example.

\section{Rationale behind algorithm}
\label{sec:rationale}
The meshing algorithm determines a conforming mesh for $\Omega$ by
perturbing vertices in the triangulation ${\cal T}_h^{0,1,2}$.
Here we briefly discuss the rationale behind the vertex adjustments we
perform--- why we perturb vertices in a particular way and when such
perturbations are possible.

\noindent
\textbf{Mapping vertices to their closest point on
  ${\boldsymbol \partial}{\boldsymbol \Omega}$}: Each vertex of a positive
edge is mapped to its closest point on the boundary. This step
transforms ${\cal T}_h^{0,1,2}$ into a conforming mesh for
$\Omega$. The reason is an intuitive one. As we discuss below,  the only
vertices in ${\cal T}_h^{0,1,2}$ that lie outside $\Omega$ are the
vertices of positive edges. Moreover, the collection of positive edges
are the boundary edges of ${\cal T}_h^{0,1,2}$ because each of these edges
belongs to just one positively cut triangle. By snapping vertices of
positive edges onto $\partial\Omega$, every vertex in the resulting
mesh belongs to $\overline{\Omega}$ and the boundary edges in the
final mesh interpolate $\partial\Omega$.

A small mesh size near $\partial\Omega$ is essential for this
perturbation step. For otherwise, $\pi$ may be multi-valued at
vertices of positive edges. A small mesh size is also required to know
that positive edges are boundary edges of ${\cal T}_h^{0,1,2}$. If the
mesh size is too large, it is possible for two
positively cut triangles to share a common positive edge. The acute
conditioning angle assumption in \S\ref{subsec:assumptions} is
critical as well--- it ensures that both vertices of a positive edge
never map to the same point in $\partial\Omega$ hence preventing
positively cut triangles from being mapped to degenerate ones. In
fact, using the assumptions in \S\ref{subsec:assumptions}, we proved
in \cite{rangarajan2011analysis} that the restriction of $\pi$ to the
collection of positive edges is a homeomorphism onto $\partial\Omega$
with Jacobian bounded away from zero. Hence we know that moving
vertices of positive edges onto the boundary does not yield a tangled
mesh. The bound for the Jacobian also implies that positive edges are
mapped to interpolating edges whose lengths are neither too small nor
too large. We refer to \cite{rangarajan2011analysis,
  rangarajan2011parameterization} for a detailed discussion on the
sufficiency of the assumptions in \S\ref{subsec:assumptions} to show
that $\partial\Omega$ is parameterized over the collection of positive
edges. There we also mention a way of relaxing the acute conditioning
angle assumption.

\begin{figure}
  \centering
  \includegraphics[scale=0.8]{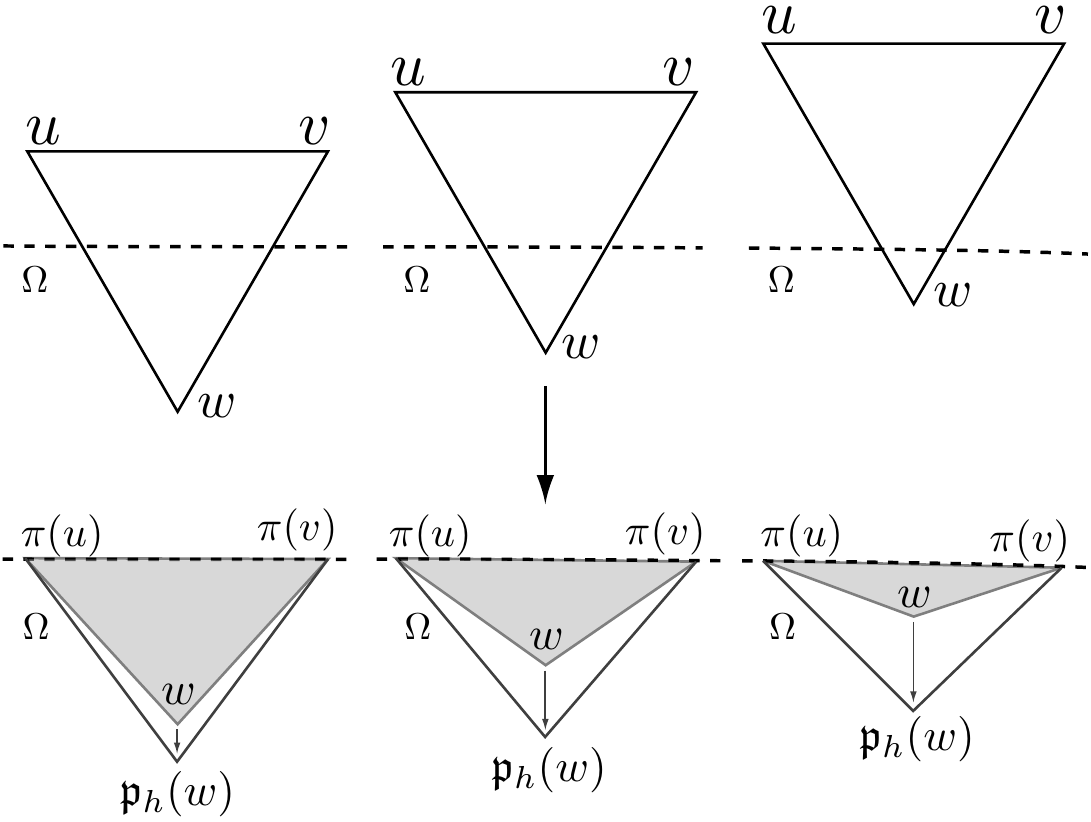}
  \caption{Vertices in $\Omega$ near its boundary need to be relaxed
    away from the boundary in the meshing algorithm to ensure good
    quality triangles. The figure shows a positively cut triangle with
    vertices $\{u,v,w\}$ and positive edge $\overline{uv}$.  Step 3
    in Table \ref{algo} moves vertices $u$ and $v$ to $\pi(u)$ and
    $\pi(v)$ respectively. The quality of the resulting triangle with
    vertices $\{\pi(u),\pi(v),w\}$(shaded in gray) depends on the
    distance of $w$ from $\partial\Omega$. The three cases in the
    figure show that the closer $w$ is to $\partial\Omega$, the poorer
    the quality of the shaded triangle. Perturbing $w$ away from
    $\partial\Omega$ alleviates this.}
  \label{fig:fail-2}
\end{figure}
\textbf{Relaxing vertices away from
  ${\boldsymbol \partial}{\boldsymbol \Omega}$}: Next we explain using
an example why it is necessary to relax vertices away from the
boundary (step 4 in Table \ref{algo}). Consider a positively cut
triangle $K$ with vertices $\{u,v,w\}$ and positive edge
$\overline{uv}$ as shown in Fig. \ref{fig:fail-2}. By snapping vertices
$u$ and $v$ onto $\partial\Omega$ in step 3 in the algorithm, $K$ is
mapped to a triangle $\tilde{K}$ with vertices $\{\pi(u),\pi(v),w\}$
(shaded in gray in the figure). We know from the above discussion that
the length of the edge $\overline{\pi(u)\pi(v)}$ cannot be too
small. However, lengths edges $\overline{\pi(u)w}$ and
$\overline{\pi(v)w}$ can be arbitrarily small.  As depicted in the
figure, the closer $w$ is to $\partial\Omega$, the poorer the quality
of triangle $\tilde{K}$. To alleviate this, we move vertex $w$ away
from $\partial\Omega$ by a small distance. In turn, to accommodate such
perturbations, vertices in a small neighborhood of $\partial\Omega$
are relaxed away from it.

In the map $\mathfrak{p}_h$ used to relax vertices in the algorithm,
the parameter $r$ determines the neighborhood of $\partial\Omega$ in
which vertices are perturbed. Since we pick $r$ to be a few multiples
of the mesh size near $\partial\Omega$, $\mathfrak{p}_h$ is only a
local perturbation near $\partial\Omega$. To prevent this step from
inducing overlapping or poorly shaped triangles, a small-enough mesh
size is also needed. For example, $\pi$ needs to be well defined in
the $r$-neighborhood of $\partial\Omega$, so that $\mathfrak{p}_h$ is
well defined as well over the vertices that are relaxed. Although a
detailed analysis is still needed, this step does not appear to pose
more stringent requirements on the mesh size that those posed by
$\varphi^h_K$ or $\Psi^h_K$.


\section{Concluding remarks}
\label{sec:conclusion}
The meshing algorithm, the mappings to curvilinear triangles used in
constructing high-order curved finite elements, and the idea of
universal meshes are useful tools for an important class of
computationally challenging problems. By employing them in problems
with evolving fluid domains while using a single background mesh, we
demonstrated the algorithmic advantages they offer.  We envision their
application to more demanding problems ranging from dynamic crack
propagation to phase transformations.

For these tools to also be useful in realistic engineering
applications, important questions remain. A significant one is knowing
what is a sufficiently small mesh size for the background mesh.  A
computable estimate is valuable because it can help determine if and
when a background(universal) mesh needs to be changed during the
course of simulating an evolving domain. The mesh size estimates in
\cite{rangarajan2011analysis} for parameterizing the immersed boundary
will be useful in determining such bounds.

A second challenge is summarized by the fact that while we can mesh
extremely complex smooth domains using a simple background mesh, we
have not specified how to mesh a square. This is essentially a
consequence of choosing the closest point projection to parameterize
the boundary. Whether domains with corners, cracks, and interfaces can
be handled without introducing additional restrictions on the
background mesh remains to be seen. An important step towards meshing
such domains is parameterizing immersed curves with end points and
corners. We have shown how to do this in
\cite{rangarajan2011parameterization}.

We think the ideas introduced here can be extended to meshing three
dimensional domains immersed in background meshes of tetrahedral
meshes. An analysis will reveal the necessary requirements on the
background mesh.

Finally, we mention that maintaining the regularity of evolving
boundaries has been a recurring challenge in numerical methods for
moving boundary problems. It requires careful choices for representing
the domain and for schemes to advance the boundary.  The literature on
these topics is growing and will continue to benefit from each new
contribution.

\FloatBarrier

\bibliographystyle{siam} 
\bibliography{references}

\begin{thebibliography}{10}

\bibitem{angot1999penalization}
{\sc P.~Angot, C.H. Bruneau, and P.~Fabrie}, {\em A penalization method to take
  into account obstacles in incompressible viscous flows}, Numerische
  Mathematik, 81 (1999), pp.~497--520.

\bibitem{babuska1990plate}
{\sc I.~Babu{\v{s}}ka and J.~Pitk{\"{a}}ranta}, {\em The plate paradox for hard
  and soft simple support}, SIAM Journal on Mathematical Analysis, 21 (1990),
  p.~551.

\bibitem{bazilevs2008isogeometric}
{\sc Y.~Bazilevs, V.M. Calo, T.J.R. Hughes, and Y.~Zhang}, {\em Isogeometric
  fluid-structure interaction: theory, algorithms, and computations},
  Computational Mechanics, 43 (2008), pp.~3--37.

\bibitem{bern1994provably}
{\sc M.~Bern, D.~Eppstein, and J.~Gilbert}, {\em Provably good mesh
  generation}, Journal of Computer and System Sciences, 48 (1994),
  pp.~384--409.

\bibitem{borgers1990triangulation}
{\sc C.~B{\"o}rgers}, {\em A triangulation algorithm for fast elliptic solvers
  based on domain imbedding}, SIAM Journal on Numerical Analysis,  (1990),
  pp.~1187--1196.

\bibitem{burman2010fictitious}
{\sc E.~Burman and P.~Hansbo}, {\em Fictitious domain finite element methods
  using cut elements: I. a stabilized lagrange multiplier method}, Computer
  Methods in Applied Mechanics and Engineering, 199 (2010), pp.~2680--2686.

\bibitem{chapelle2003finite}
{\sc D.~Chapelle and K.J. Bathe}, {\em The finite element analysis of shells:
  fundamentals}, Springer Verlag, 2003.

\bibitem{ciarlet1978finite}
{\sc P.G. Ciarlet}, {\em The finite element method for elliptic problems},
  vol.~4, North Holland, 1978.

\bibitem{ciarlet1972interpolation}
{\sc P.G. Ciarlet and P.A. Raviart}, {\em Interpolation theory over curved
  elements, with applications to finite element methods}, Computer Methods in
  Applied Mechanics and Engineering, 1 (1972), pp.~217--249.

\bibitem{donea2004arbitrary}
{\sc J.~Donea, A.~Huerta, J.P. Ponthot, and A.~Rodr{\'\i}guez-Ferran}, {\em
  Arbitrary {Lagrangian}--{Eulerian} methods}, Encyclopedia of Computational
  Mechanics,  (2004).

\bibitem{girault1986finite}
{\sc V.~Girault and P.~Raviart}, {\em Finite element methods for Navier-Stokes
  equations: Theory and Algorithms}, vol.~5 of Springer Series in Computational
  Mathematics, Springer-Verlag, Berlin and New York, 1986.

\bibitem{gonzalez2006inverse}
{\sc M.~Gonzalez and M.~Goldschmit}, {\em Inverse geometry heat transfer
  problem based on a radial basis functions geometry representation},
  International journal for numerical methods in engineering, 65 (2006),
  pp.~1243--1268.

\bibitem{gordon1973transfinite}
{\sc W.J. Gordon and C.A. Hall}, {\em Transfinite element methods:
  blending-function interpolation over arbitrary curved element domains},
  Numerische Mathematik, 21 (1973), pp.~109--129.

\bibitem{hansbo2002unfitted}
{\sc A.~Hansbo and P.~Hansbo}, {\em An unfitted finite element method, based on
  nitsche's method, for elliptic interface problems}, Computer methods in
  applied mechanics and engineering, 191 (2002), pp.~5537--5552.

\bibitem{lenoir1986optimal}
{\sc M.~Lenoir}, {\em Optimal isoparametric finite elements and error estimates
  for domains involving curved boundaries}, SIAM Journal on Numerical Analysis,
   (1986), pp.~562--580.

\bibitem{lew2008discontinuous}
{\sc A.J. Lew and G.C. Buscaglia}, {\em A discontinuous-{Galerkin}-based
  immersed boundary method}, International Journal for Numerical Methods in
  Engineering, 76 (2008), pp.~427--454.

\bibitem{mansfield1978approximation}
{\sc L.~Mansfield}, {\em Approximation of the boundary in the finite element
  solution of fourth order problems}, SIAM Journal on Numerical Analysis,
  (1978), pp.~568--579.

\bibitem{moumnassi2011finite}
{\sc M.~Moumnassi, S.~Belouettar, E.~Bechet, S.~Bordas, D.~Quoirin, and
  M.~Potier-Ferry}, {\em Finite element analysis on implicitly defined domains:
  An accurate representation based on arbitrary parametric surfaces}, Computer
  methods in applied mechanics and engineering, 200 (2011), pp.~774--796.

\bibitem{osher2003level}
{\sc S.~Osher and R.P. Fedkiw}, {\em Level set methods and dynamic implicit
  surfaces}, vol.~153, Springer Verlag, 2003.

\bibitem{peskin2002immersed}
{\sc C.S. Peskin}, {\em The immersed boundary method}, Acta Numerica, 11
  (2002), pp.~479--517.

\bibitem{rangarajan2011analysis}
{\sc R.~Rangarajan and A.J. Lew}, {\em Analysis of a method to parameterize
  planar curves immersed in triangulations}, ArXiv e-prints,  (2011).

\bibitem{rangarajan2009discontinuous}
{\sc R.~Rangarajan, A.J. Lew, and G.C. Buscaglia}, {\em A
  discontinuous-{Galerkin}-based immersed boundary method with non-homogeneous
  boundary conditions and its application to elasticity}, Computer Methods in
  Applied Mechanics and Engineering, 198 (2009), pp.~1513--1534.

\bibitem{rangarajan2011parameterization}
{\sc Ramsharan Rangarajan and Adrian~J. Lew}, {\em Parameterization of planar
  curves immersed in triangulations with application to finite elements},
  International Journal for Numerical Methods in Engineering, 88 (2011),
  pp.~556--585.

\bibitem{saksono2007adaptive}
{\sc P.~Saksono, W.~Dettmer, and D.~Peri{\'c}}, {\em An adaptive remeshing
  strategy for flows with moving boundaries and fluid--structure interaction},
  International Journal for Numerical Methods in Engineering, 71 (2007),
  pp.~1009--1050.

\bibitem{sanches2011immersed}
{\sc R.~Sanches, P.~Bornemann, and F.~Cirak}, {\em Immersed b-spline (i-spline)
  finite element method for geometrically complex domains}, Computer Methods in
  Applied Mechanics and Engineering,  (2011).

\bibitem{scott1973finite}
{\sc R.~Scott}, {\em Finite element techniques for curved boundaries}, PhD
  thesis, Massachusetts Institute of Technology, Cambridge, 1973.

\bibitem{wagner2001extended}
{\sc G.~Wagner, N.~Mo{\"e}s, W.~Liu, and T.~Belytschko}, {\em The extended
  finite element method for rigid particles in stokes flow}, International
  Journal for Numerical Methods in Engineering, 51 (2001), pp.~293--313.

\bibitem{zlamal1973curved}
{\sc M.~Zl{\'a}mal}, {\em Curved elements in the finite element method. {I}},
  SIAM Journal on Numerical Analysis,  (1973), pp.~229--240.

\end{thebibliography}

\appendix
\section{Implementation for the meshing algorithm}
\label{sec:algo}
We provide a simple implementation of the meshing algorithm in
\S\ref{sec:meshing} to determine a conforming mesh for $\Omega$.  We
assume that the background mesh ${\cal T}_h$ is specified by (i) a
list of coordinates $V$ of its vertices, (ii) a numbering $I$ for the
vertices, and (iii) a list of triangle connectivities ${C}$ which are
$3$-tuples of vertex numbers. Vertex with number $i\in {I}$ is denoted
$v_i\in V$. A triangle with vertices $\{v_i,v_j,v_k\}$ has
connectivity $(i,j,k)\in C$. The algorithm returns a conforming
triangulation for $\Omega$ specified by a set of vertices $V^\Omega$,
a numbering $I^\Omega\subseteq {I}$ for these vertices and a
connectivity list ${C}^\Omega\subseteq C$ for triangles in the mesh.

\algsetup{indent=2em}
\begin{algorithmic}[1]
  \REQUIRE vertex coordinates $V$, vertex numbering ${I}$, triangle
  connectivities ${C}$.

  \REQUIRE Choose $\eta\in (0,1)$ and $r$ equal to a few multiples of
  $h$
    
  \STATE Initialize $I^\Omega \leftarrow \emptyset, V^\Omega\leftarrow
  \emptyset, C^\Omega\leftarrow \emptyset$
    
    \FORALL{$i\in I$}
    \IF{$v_i\in\Omega$}
    \STATE $s_i\leftarrow \text{true}$
    \ELSE
    \STATE $s_i\leftarrow \text{false}$
    \ENDIF
    \IF{$s_i$}
    \STATE Append $i$ to $I^\Omega$
    \IF{$-r < \phi(v_i)$}
    \STATE Append $\mathfrak{p}_h(v_i)$ to $V^\Omega$
    \ELSE
    \STATE Append $v_i$ to $V^\Omega$
    \ENDIF
    \ENDIF
    \ENDFOR
    
    \STATE $I^+\leftarrow \emptyset$
    \FORALL{$(i,j,k)\in C$}
    \STATE $\mathfrak{i}_- \leftarrow \emptyset, \mathfrak{i}_+\leftarrow \emptyset$
    \FOR{$\ell\in \{i,j,k\}$}
    \IF{$s_\ell$}
    \STATE Append $\ell$ to  $\mathfrak{i}_-$
    \ELSE 
    \STATE Append $\ell$ to $\mathfrak{i}_+$
    \ENDIF
    \ENDFOR
    
    \IF{$\#\mathfrak{i}_-\geq 1$}
    \STATE Append $(i,j,k)$ to $C^\Omega$
    \IF{$\#\mathfrak{i}_+=2$}
    \STATE Ensure conditioning angle of triangle $(i,j,k)$ is acute.
    \FOR{$\ell\in \mathfrak{i}_+$}
    \STATE Append $\ell$ to $I^+$
    \ENDFOR
    \ENDIF
    \ENDIF
    
    \ENDFOR

    \FORALL{$i\in I^+$}
    \STATE Append $i$ to $I^\Omega$
    \STATE Append $\pi(v_i)$ to $V^\Omega$
    \ENDFOR
    
    \RETURN triangulation $(V^\Omega, I^\Omega, C^\Omega)$ for $\Omega$
\end{algorithmic}

\end{document}